\documentclass[11pt]{article}
\usepackage{amsmath}
\usepackage{graphicx}
\usepackage{epsfig}
\usepackage{amsfonts}
\usepackage{amssymb}
\usepackage{placeins}
\usepackage{cases}
\usepackage[margin=1in]{geometry}
\usepackage{authblk}
\usepackage[titletoc,toc,title]{appendix}
\usepackage{epstopdf}
\usepackage{changes}
\usepackage{xcolor}
\usepackage{bm}
\usepackage{enumerate}

\usepackage{accents}

\title{Reduced Order Characterization of Nonlinear Oscillations Using an Adaptive Phase-Amplitude Coordinate Framework}

\begin{document}
\author[1]{Dan Wilson \thanks{corresponding author:~dwilso81@utk.edu}}
\author[1]{Kai Sun}
\affil[1]{Department of Electrical Engineering and Computer Science, University of Tennessee, Knoxville, TN 37996, USA}
\maketitle


\begin{abstract}
We propose a general strategy for reduced order modeling of systems that display highly nonlinear oscillations.  By considering a continuous family of forced periodic orbits defined in relation to a stable fixed point and subsequently leveraging phase-amplitude-based reduction strategies, we arrive at a low order model capable of accurately capturing nonlinear oscillations resulting from arbitrary external inputs.  In the limit that oscillations are small, the system dynamics relax to those obtained from local linearization, i.e.,~that can be fully described using linear eigenmodes.  For larger amplitude oscillations, the behavior can be understood in terms of the dynamics of a small number of nonlinear modes.  We illustrate the proposed strategy in a variety of examples yielding results that are substantially better than those obtained using standard linearization-based techniques.

\end{abstract}

\section{Introduction}
We develop and investigate a general strategy for reduced order representation of systems displaying highly nonlinear oscillations.  This is in direct contrast to linear oscillations which can be decomposed into a superposition of eigenmodes with growth/decay/oscillation rates governed by their associated eigenvalues \cite{ewin09}, \cite{hesp09}.  While such approaches can be used to approximate low amplitude oscillations for a nonlinear system, alternative techniques must be used to accommodate large amplitude oscillations occurring in nonlinear systems.

Early work on the subject of nonlinear oscillations considered the idea of a nonlinear normal mode \cite{rose60}, \cite{rose66} as a synchronous oscillation admitted by a nonlinear system which provided a nonlinear extension to linear modes.  Subsequent work  \cite{shaw91}, \cite{shaw94} viewed these nonlinear normal modes as two-dimensional invariant manifolds that are tangent to a related linear eigenspace.  Review articles \cite{kers09} and \cite{avra10} discuss more recent results and generalizations of the use of nonlinear normal modes in the analysis of nonlinear oscillations.  Further extensions were considered in \cite{hall16} which introduced the notion of a spectral submanifold, the smoothest invariant manifold that functions as a nonlinear extension of a linear modal subspace.  Investigation of the dynamics on these spectral submanifolds can yield information about nonlinear oscillations and extract associated backbone curves \cite{pons20}, \cite{szal17}.

From a broader perspective, reduced order modeling of nonlinear oscillations can be thought of as a dynamical representation problem \cite{mezi21}.   Koopman-based approaches, for instance, attempt to represent the dynamics of a general nonlinear system using a linear, but possibly infinite dimensional operator \cite{budi12}, \cite{mezi13}, \cite{mezi19}.  In direct contrast to local linearization techniques that consider the dynamics in a close neighborhood of some nominal solution, Koopman-based approaches can be used to obtain linear representations for the fully nonlinear dynamics of observables, subsequently allowing for the analysis of nonlinear oscillations in terms of the superposition Koopman eigenmodes.  From a practical perspective, the key challenge of implementing Koopman-based approaches is in the identification of a suitable finite basis to represent the possibly infinite dimensional Koopman operator.  In some cases, this can be accomplished by finding a Koopman invariant subspace to yield an exact linear, finite-dimensional representation for a nonlinear system \cite{take17}, \cite{kord20}, \cite{brun16}.  More commonly, data-driven algorithms such as dynamic mode decomposition \cite{schm10}, \cite{kutz16}, \cite{will15} are used to provide finite dimensional, linear approximations for the action of the Koopman operator.

Rather than characterizing the full action of the Koopman operator, a number of authors have suggested the use of a subset of Koopman eigenfunctions to establish a reduced order coordinate system for representing the dynamics of a fully nonlinear system \cite{maur13}, \cite{wils16isos}, \cite{kais21}, \cite{kval21}. Among these is the isostable coordinate framework \cite{maur13}, which considers the level sets of the slowest decaying Koopman eigenmodes to form a reduced order basis.  Previous work \cite{wils20highacc}, \cite{wils21dd}, \cite{wils20ddred}, \cite{wils17isored}  has considered this general coordinate system in the development of various model order reduction algorithms that are applicable to systems with fixed points and periodic orbits.  By retaining only the slow decaying components and truncating the rest, high accuracy reduced order models can be obtained that contain only a handful of state variables. 

In conjunction with the isostable coordinate basis, recent work \cite{wils21adapt}, \cite{wils21adaptpde} proposed an adaptive coordinate system that considers a family of either stable limit cycles or fixed points that emerge when using different parameter sets.  By adaptively selecting the nominal attractor from within this family, provided the state remains close to the  attractor (i.e.,~as gauged by the magnitude of the isostable coordinates) truncation errors can be mitigated resulting in a very accurate but still substantially reduced order models.  While initial results obtained using this  adaptive reduction strategy have been promising \cite{wils21optimal}, \cite{wils22adapt}, \cite{toth22}, many unanswered questions remain regarding its implementation.  For instance, there is usually no systematic way of choosing the family of reference trajectories.  Additionally, while some heuristics are discussed in \cite{wils21optimal} for adaptively selecting the nominal attractor in order to limit truncation errors, it is not always obvious how to accomplish this task.

Here, we consider a general strategy for characterizing nonlinear oscillations for systems with stable fixed points using the aforementioned adaptive phase-amplitude reduction approach.  As a primary contribution, this work proposes and investigates a systematic strategy for defining reference trajectories used in conjunction with the adaptive phase-amplitude reduction strategy.  The resulting approach yields a reduced order model that can consider large amplitude oscillations that result from arbitrary forcing far beyond the linear regime.  In the limit that oscillations are small, the system dynamics relax to those obtained from local linearization.  Additionally, the proposed approach can also be used to explicitly consider the interactions between multiple nonlinear modes.   The organization of this paper is as follows:~Section \ref{backsec} provides necessary background on phase and phase-amplitude reduction techniques that are used as a starting point for the proposed model order reduction strategy.  Section \ref{derivsec} provides a description of the proposed approach and provides a detailed list of steps required for implementation.  Section \ref{ressec} provides numerical illustrations for three example systems:~a simple pendulum, a coupled population of nonlinear planar oscillators, and a power system model comprised of synchronous generators. In each example, the proposed strategy far exceeds others that leverage techniques based on local linearization.  Section \ref{concsec} provides concluding remarks.

\section{Background on Phase-Based Reduction Techniques for Oscillatory Systems} \label{backsec}

Consider an ordinary differential equation of the form
\begin{equation} \label{maineq}
\dot{x} = F(x,u),
\end{equation}
where $F$ sets the generally nonlinear dynamics, $x \in \mathbb{R}^N$ is the state, and $u \in \mathbb{R}^M$ is an input.  Suppose that when $u = 0$, Equation \eqref{maineq} has a stable fixed point $x_0$ for which $F(x_0,0) = 0$.  Letting $\Delta x = x - x_0$, a standard linearization yields a model of the form 
\begin{equation}\label{lineq}
\Delta \dot{x} = A \Delta x + B u,
\end{equation}
where $A = \frac{\partial F}{\partial x}$ and $B = \frac{\partial F}{\partial u}$, both evaluated at $x = x_0$ and $u = 0$.  Assume that all eigenvalues of $A$ are simple, i.e.,~unique.  In a close neighborhood of the fixed point, when $u = 0$, Equation \eqref{lineq} admits solutions of the form
\begin{equation} \label{evaldecomp}
\Delta x (t) = \sum_{j = 1}^N  s_j v_j \exp(\lambda_j t) + O(||\Delta x||^2),
\end{equation}
where $(\lambda_j,v_j)$ is an eigenvalue/eigenvector pair of the matrix $A$ and $s_j$ is the associated coefficient after projecting $\Delta x$ onto the eigenspace.  The decomposition \eqref{evaldecomp} can be used to analyze  linear oscillations in \eqref{lineq} in terms of the eigenvectors associated with complex-valued eigenvalues.  The influence of nonzero input can also be readily considered with a variety of linear control techniques \cite{skog05}, \cite{hesp09}.  However, this linearization is only valid in the limit that $\Delta x$ and $u$ are both small, precluding use in many practical applications.  In order to overcome this limitation, the strategy proposed in this work considers the application of additional periodic forcing and the resulting periodic orbits to characterize the dynamics far from the fixed point.  This will be done using phase-amplitude-based reduced order modeling approaches as a starting point which are briefly summarized below.

\subsection{Phase Reduction} 
Consider a general ordinary differential equation of the form
\begin{equation} \label{phaseeq}
\dot{y} = F(y,p) + u(t),
\end{equation}
where $y \in \mathbb{R}^N$ is the state, $F$ gives the dynamics, $p \in \mathbb{R}^M$ is a parameter set, and $u$ is an additive input.  Suppose that for a constant choice of $p = q$ and when $u = 0$, Equation \eqref{phaseeq} has a stable $T(q)$-periodic orbit $y_q^\gamma$.  Isochrons \cite{guck75}, \cite{winf01} can be used to define an oscillation phase at all locations in the basin of attraction of this limit cycle.  Letting $\theta_1 \in [0,2\pi)$ be the phase corresponding to some initial condition $a(0) \in y^\gamma_q$, the $\theta_1$ isochron is given by the set of all $b(0)$ for which
\begin{equation} \label{isodef}
\lim_{t \rightarrow \infty} || a(t)-b(t)|| = 0,
\end{equation}
where $||\cdot||$ is some vector norm.  Considering Equation \eqref{isodef}, initial conditions that share the same isochron share the same asymptotic convergence to the periodic orbit.  Typically, the phase $\theta(y,q)$ is scaled so that $\frac{d \theta}{dt} = \frac{2 \pi}{T(q)} = \omega(q)$ under the flow of \eqref{phaseeq}.  With the definition of phase in mind, a standard phase reduction \cite{erme10}, \cite{kura84}, \cite{winf01}, can be performed by restricting attention to a close neighborhood of the periodic orbit and changing variables to phase coordinates 
\begin{align} \label{predeq}
\frac{d \theta}{dt} &= \frac{\partial \theta}{\partial y} \cdot \frac{dy}{dt} \nonumber \\
&=  \frac{\partial \theta}{\partial y} \cdot \Big(  F(y,q) + u(t) \Big) \nonumber \\
&= \omega(p) + Z^T(\theta,q) u(t) + O(||y - y^\gamma_q(\theta)||^2),
\end{align}
where $Z(\theta,q) = \frac{\partial \theta}{\partial y}$ evaluated at $y^\gamma_q(\theta)$, the dot denotes the dot product, and $^T$ denotes the transpose.  Above, $\frac{\partial \theta}{\partial x} \cdot F(y,q) = \omega$ because $\frac{d \theta}{dt} = \omega$ when $u = 0$.  Typically, the higher order terms are truncated to yield a closed form ordinary differential equation that can be used to represent the behavior of the original $N$-dimensional system in terms of a 1-dimensional reduction.  Phase reductions of the form \eqref{predeq} have been used extensively to analyze and understand dynamics that emerge in weakly perturbed oscillatory systems \cite{erme10}, \cite{brow03}, \cite{wils14a}, \cite{schw12}, \cite{piet19}.

\subsection{Phase-Amplitude Reduction}
Equation \eqref{predeq} is only valid in the limit of weak forcing.  When considering large magnitude inputs, information about amplitude coordinates must be considered which represent directions transverse to the periodic orbit.  There are wide a variety of strategies that can be used to incorporate the influence of amplitude-based effects \cite{wils16isos}, \cite{wedg13}, \cite{lets20}, \cite{wils18operat}, \cite{kure13}, \cite{cast13}, \cite{rose19}.  Here, we will use Floquet coordinates for this task.  Considering the periodic orbit of \eqref{phaseeq} and letting $\Delta y = y- y^\gamma_q(\theta)$, to a linear approximation one can write $\Delta \dot{y} = J \Delta y$ where $J$ is the Jacobian evaluated at $x^\gamma_q(\theta(t))$.   Defining $\Phi$ to be the monodromy matrix of this $T(p)$-periodic linear time varying system, provided $\Phi$ is diagonalizable, near $y^\gamma_q$ one can leverage Floquet theory \cite{jord07} to write
\begin{equation}  \label{floqexp}
y - y^\gamma_q(\theta) = \sum_{j = 1}^{N-1} \psi_j g_j(\theta,q) + {O}(\psi_1^2) + \dots + {O}(\psi_{N-1}^2),
\end{equation} 
where $g_j(\theta,q)$ is a Floquet eigenfunction and $\psi_1,\dots,\psi_{N-1}$ are associated Floquet coordinates.  Note that in \eqref{floqexp}, the contribution from the $N$th Floquet eigenfunction is absorbed by the phase coordinate.  Equation \eqref{floqexp} can also be extended to nonlinear orders of accuracy using the notion of isostable coordinates \cite{wils20highacc}, \cite{wils17isored} which can be defined in the entire basin of attraction of the limit cycle.  To linear order, the Floquet eigenfunctions from Equation \eqref{floqexp} can be used to augment the phase reduction from \eqref{predeq} yielding a phase-amplitude reduction
\begin{align} \label{phaseamp}
\dot{\theta} &= \omega(q) + Z^T(\theta,q) u(t), \nonumber \\
\dot{\psi}_j &= \kappa_j(q) \psi_j + I_j^T(\theta,q) u(t),  \nonumber  \\
j & = 1,\dots,N-1.
\end{align}
Here $I_j(\theta,q) = \frac{\partial \psi_j}{\partial y}$ evaluated at $y^\gamma_q(\theta)$ and $\kappa_j$ is the Floquet exponent corresponding to the $j^{\rm th}$ Floquet eigenfunction.  While Equation \eqref{phaseamp} is still only valid provided the state $y$ of the underlying Equation \eqref{phaseeq} remains close to the underlying periodic orbit $y^\gamma_q$.  Nonetheless, this additional information can still be useful, for instance, in the context of control design \cite{mong19b}, \cite{wils16isos}.  Note that references such as \cite{wils20highacc} and \cite{wils17isored} refer to $\psi_1,\dots,\psi_{N-1}$ as isostable coordinates; because isostable coordinates are identical to Floquet coordinates to linear order, and because this work only considers the dynamics of these coordinates to linear order for a given periodic orbit, we will refer to these amplitude coordinates as Floquet coordinates in this manuscript.

\subsection{Adaptive Phase-Amplitude Reduction} \label{predsec}

 If it is necessary to consider inputs that drive the system far beyond its reference periodic orbit, the adaptive phase-amplitude reduction can be used  \cite{wils21adapt}, \cite{wils21optimal}.  To implement this strategy, assume that for some allowable $q \in \mathbb{R}^M$, when $q$ is held constant the equation $\dot{y} = F(y,q)$ has a periodic orbit $y^\gamma_q$.  For each of these orbits, one can define an extended phase $\theta(x,q)$ and a set of extended Floquet coordinates $\psi_1(x,q),\dots,\psi_{N-1}(x,q)$.  Note that the phase coordinates are unique to a constant shift and must be disambiguated, for instance, by choosing the crossing of some Poincar\'e section to correspond to a level set of phase for each limit cycle.   Equation \eqref{phaseeq} can subsequently be rewritten as
\begin{equation} \label{yeq}
\dot{y} = F(y,q) + U_e(y,u,p,q),
\end{equation}
where 
\begin{equation}
U_e(y,u,p,q) = F(y,p) - F(y,q) + u(t).
\end{equation}
Rewriting in this manner, the dynamics can be considered relative to the reference orbit $y^\gamma_q$ with effective input $U_e$.  As described in \cite{wils21adapt}, allowing $q$ to be nonstatic (i.e.,~allowing the reference orbit to change), transforming to phase and Floquet coordinates and truncating  all $O(||y - y^\gamma_q(\theta)||^2)$ terms yields
\begin{align} \label{phaseampred}
\dot{\theta} &= \omega(q) + Z^T(\theta,q) U_e + D^T(\theta,q) \dot{q}, \nonumber \\
\dot{\psi}_j &= \kappa_j(q) \psi_j + I_j^T(\theta,q)U_e + E^T_j(\theta,q) \dot{q}, \nonumber \\
j & =1, \dots,N-1, \nonumber \\
\dot{q} &= G_q(q,\theta,\psi_1,\dots,\psi_{N-1},U_e).
\end{align}
Above, $D(\theta,q) \in \mathbb{R}^M$ with the $i^{\rm th}$ element capturing how changes to $q$ yield changes in $\theta$, similarly, each $E_j(\theta,q) \in \mathbb{C}^M$ characterizes how changes to $q$ influence the Floquet coordinates.  More details about the computation of these terms are provided in \cite{wils21adapt}.  The term $G_q$ sets the update rule for the parameter set $q$.  Equation \eqref{phaseampred} is valid in the limit that $||y - y^\gamma_q(\theta)||$ remains small.  Considering Equation \eqref{floqexp}, provided that the Floquet coordinates $\psi_1,\dots,\psi_{N-1}$ can be kept small, $y - y^\gamma_q$ will remain small.  As such, by designing $G_q$ such that each Floquet coordinate remains small, the adaptive reduction \eqref{phaseampred} can be used to accurately represent dynamics of the underlying system \eqref{phaseeq}, even when the inputs considered are large.  

Note that Equation \eqref{phaseampred} is actually higher dimensional than the original equation \eqref{phaseeq} ($N+M$ dimensions versus $N$ dimensions).  To obtain a reduced order equation, in many cases, a large number of Floquet coordinates can be well approximated by zero if their corresponding Floquet exponents are negative and large in magnitude so that they decay rapidly \cite{wils16isos}, \cite{wils21adapt}.  In other cases, a coordinate $\psi_j$ can be neglected when the input $U_e$ is orthogonal to $I_j(\theta,p)$ \cite{wils19phase}.  General heuristics for choosing $G_p$ are discussed in \cite{wils21adapt}.

\section{Derivation of the Proposed Technique for Defining Nonlinear Oscillatory Modes}  \label{derivsec}

\subsection{Overview of the Proposed Strategy} \label{overviewsec}
The overall goal is to characterize oscillations in the model \eqref{maineq} in terms of an appropriate subset of nonlinear modes that are valid far beyond the weakly perturbed limit.  Towards a formulation that leverages the adaptive phase-amplitude reduction described in Section \ref{predsec}, consider the same system with an additional input
\begin{equation} \label{fulleq}
\dot{x} = F(x,u) + \alpha(q,t).
\end{equation}
Here, $x$, $F$, and $u$ represent the state, dynamics, and input as defined in Equation \eqref{maineq}, $q \in \mathbb{R}^M$ is a time-varying parameter set, and $\alpha(q,t)$ is an external periodic input with period $T$.   Suppose that when $u = \alpha = 0$, a fixed point $x_0$ exists.   Letting $s \in [0,T)$, Equation \eqref{fulleq} can be rewritten as an autonomous system of equations:
\begin{align} \label{nonaut}
\dot{x} &= F(x,u) + \alpha(q,s), \nonumber \\
\dot{s} &= 1.
\end{align}
Letting $y = [x^T \; s]^T \in \mathbb{R}^{N+1}$, suppose that when $u = 0$, for all allowable values of $q$ Equation \eqref{nonaut} has a periodic orbit
\begin{equation} \label{yorb}
y^\gamma_q(t) = \begin{bmatrix} x^\gamma_q(t) \\ {\rm mod}(t,T) \end{bmatrix}.
\end{equation}
Following the formulation from Section \ref{predsec}, for each $y^\gamma_q$, an extended phase $\theta(x,q)$ and set of Floquet coordinates $\psi_1(x,q),\dots,\psi_{N}(x,q)$ can be defined (note that there are $N$ Floquet coordinates because $y \in \mathbb{R}^{N+1}$).  One can consider a rewritten version of Equation \eqref{nonaut}
\begin{align}  \label{adaptiveeq}
\dot{x} &= F(x,0) + \alpha(q,s) + U_e(x,q,s), \nonumber \\
\dot{s} &= 1 + f_\theta(x,q,s),
\end{align}
where $f_\theta(x,q,s) \in \mathbb{R}$ is a function with form that will be discussed momentarily and 
\begin{equation} \label{uedef}
U_e(x,q,s) = F(x,u) - F(x,0) - \alpha(q,s).
\end{equation}
Above, the $x$ dynamics in Equation \eqref{adaptiveeq} are identical those that govern \eqref{maineq}.  Noticing that Equation \eqref{adaptiveeq} can be written in the same general form as Equation \eqref{yeq}, first considering the phase coordinates, one can write 
\begin{equation} \label{phaseadapt}
\dot{\theta} = \omega(q) + Z_1^T(\theta,q)U_e + Z_2(\theta,q) f_\theta + D^T(\theta,q) \dot{q},
\end{equation}
where $Z_1(\theta,q) \in \mathbb{R}^{N}$ and $Z_2(\theta,q) \in \mathbb{R}^{1}$ comprise the first $N$ elements and last element, respectively, of $Z(\theta,q)$  associated with the periodic orbit $y^\gamma_q(t)$ and $\omega = 2\pi/T$.  To further simplify \eqref{phaseadapt}, note that the periodic orbit from Equation \eqref{nonaut} is driven by the forcing $\alpha(q,s)$.  Because this periodic orbit emerges as the result of periodic forcing, when $q$ is held constant, one can show that $\theta(q,s) = {\rm mod}(\theta_0 + \omega s,2\pi)$ (cf.,~\cite{wils19complex}) where $\theta_0$ is an arbitrary constant.  As such, $Z_2(\theta,q) = \partial \theta / \partial s = \omega(q)$ and $Z_1(\theta,q)$ is an appropriately sized vector of zeros.  Because $\theta_0$ is arbitrary, it is most convenient to take $\theta_0$ = 0 for all $q$ yielding $D(\theta,q) = 0$.  With these simplifications in mind, Equation \eqref{phaseadapt} becomes $\dot{\theta} = \omega(q) + \omega(q) f_\theta$.  Considering the amplitude coordinates, the full adaptive phase-amplitude reduction of the form \eqref{phaseampred} is
\begin{align} \label{forcedadaptive}
\dot{\theta} &=  \omega(q) + \omega(q) f_\theta, \nonumber \\
\dot{\psi_j} &= \kappa_j(q) \psi_j + I^T_{j,1}(\theta,q)U_e + I_{j,2}(\theta,q) f_\theta + E_j^T(\theta,q) \dot{q}, \nonumber \\
& j = 1,\dots,\beta, \nonumber \\
\dot{q} &= G_q(q,\theta,\psi_1,\dots,\psi_{N},U_e).
\end{align}
Above, each $I_{j,1}(\theta,q) \in \mathbb{C}^N$ and $I_{j,2}(\theta,q) \in \mathbb{C}$ comprise the first $N$ elements and last element, respectively, of $I_j(\theta,q)$ associated with the periodic orbit $y^\gamma_q(t)$.  In \eqref{forcedadaptive}, it is assumed that the Floquet coordinates $\psi_{\beta + 1}, \dots, \psi_N$ have large magnitude Floquet exponents so that they decay rapidly and can be ignored.  

In the derivation of Equation \eqref{forcedadaptive}, the external forcing $\alpha(q,t)$ is arbitrary.   As shown in the following sections, with an appropriate choice of $\alpha(q,t)$, Equation \eqref{forcedadaptive} yields a greatly simplified, reduced order set of ordinary differential equations where the periodic  orbits $x_q^\gamma$ have a close connection with linear modes for small magnitude forcing but can be extended far beyond the linear regime.

\subsection{Construction of Appropriate Trajectories For a Single Nonlinear Oscillation Mode Near the Fixed Point} \label{singmode}

The analysis and derivations to follow in Sections \ref{singmode}-\ref{nlinmode} consider a 1-parameter family of periodic orbits, i.e.,~taking $q \in \mathbb{R}^1$, in the consideration of a single nonlinear oscillation mode.  If additional nonlinear modes are considered, an $M$-parameter family of solutions must be considered with $M>1$.  This situation is discussed in Section \ref{multimodes}.

 To begin, suppose that when $u = \alpha = 0$ in Equation \eqref{fulleq}, a fixed point $x_0$ exists for which $F(x_0,0) = 0$ with solutions given by the eigendecomposition from Equation \eqref{evaldecomp}.  Let $(\lambda_{1},v_1)$ and $(\lambda_2,v_2)$ be a simple complex-conjugate eigenvalue/eigenvector pair associated with a single oscillatory linear mode.  For simplicity, it will be assumed the eigenvalues are chosen so that ${\rm Imag}(\lambda_1) > 0$. These eigenvalues will be normalized so that $||v_1||_2 = ||v_2||_2=1$ where $||\cdot||_2$ is the 2-norm.   Note that this still leaves one additional degree of freedom since $v_1 \exp(i \phi)$ where $i = \sqrt{-1}$ is still an eigenvector for any $\phi \in \mathbb{R}$.  To fully specify $v_1$ we will require ${\rm arg}(e_j^T v_1) =-\pi$ where $e_j$ is the $j^{\rm th}$ element of the standard unit basis, ${\rm arg}(\cdot)$ is the argument of the complex number, and $j$ can be chosen arbitrarily.  Note that other normalizations are also possible. 
 
 To proceed,  let $w_1$ and $w_2$ be left eigenvectors associated with right eigenvectors $v_1$ and $v_2$.  Defining, $\mu_i(t) = w_i^T (x(t) - x_0)$ in a hyperplane orthogonal to the eigenvectors $v_3, \dots v_N$, solutions take the form
\begin{equation} \label{linmode1}
x(t) = x_0 + \sum_{j = 1}^2 \mu_j(0) v_j \exp(\lambda_j t).
\end{equation}
For simplicity, here it will be assumed that $\psi_j(0) \in \mathbb{R}$.  Considering trajectories evolving according to \eqref{linmode1}, level sets of $|\mu_1|$ are traced out by the periodic orbit
\begin{equation} \label{levelsets}
x(t) = x_0 + \sum_{j = 1}^2 \mu_j(0) v_j \exp( i [{\rm Imag}(\lambda_j) + \Delta \omega {\rm sign}({\rm Imag}(\lambda_j))] t).
\end{equation}
Above, the term $\Delta \omega$  influences the period of the orbits defined by Equation \eqref{levelsets}. It will be necessary that $\Delta \omega \neq 0$ and  $- {\rm Imag}(\lambda_1)/3 < \Delta \omega < {\rm Imag}(\lambda_1)$ so that the eigenvalues of the monodromy matrix associated with the resulting periodic orbit are simple; this point will become more clear shortly.  Next, considering Equation \eqref{fulleq}  taking $u = 0$ but $\alpha(t) \neq 0$, solutions of the form \eqref{levelsets} can be obtained by solving
\begin{align} \label{alphasolve}
 F(x_0,0) +&  \frac{\partial F}{\partial x} (x(t)-x_0) + \alpha(t) + O(||x - x_0||^2) = \frac{dx}{dt} \nonumber \\
 & = \sum_{j = 1}^2 \mu_j(0) v_j  i [{\rm Imag}(\lambda_j) + \Delta \omega {\rm sign}({\rm Imag}(\lambda_j))] \exp( i [{\rm Imag}(\lambda_j) + \Delta \omega {\rm sign}({\rm Imag}(\lambda_j))] t),
\end{align}
where the left hand side is an approximation of \eqref{fulleq} for small values of $x - x_0$ and the right hand side is obtained from direct differentiation of \eqref{levelsets}.  Solving \eqref{alphasolve} for $\alpha(t)$ and simplifying yields
\begin{align}
\alpha(t) &= -2 \mu_1(0) {\rm Real}(\lambda_1) \Big[ {\rm Real}(v_1) \cos( ({\rm Imag}(\lambda_1) + \Delta \omega) t)   \nonumber \\
& \qquad  - {\rm Imag}(v_1) \sin( ({\rm Imag}(\lambda_1) + \Delta \omega) t)  \Big]   \nonumber \\
&\qquad - 2 \Delta \omega \mu_1(0) \Big[ {\rm Imag}(v_1) \cos( ({\rm Imag}(\lambda_1) + \Delta \omega) t) \nonumber \\
&\qquad  + {\rm Real}(v_1) \sin( ({\rm Imag}(\lambda_1) + \Delta \omega) t)   \Big]     +     O(||x - x_0||^2).
\end{align}
These forced periodic orbits can be used to define a family of forced periodic orbits for use with the adaptive phase-amplitude reduction from \eqref{forcedadaptive} in the weak forcing limit.  To this end letting $q = \mu_1(0) = O(\epsilon)$ where $0 < \epsilon \ll 1$, take 
\begin{align}  \label{inporb}
 x^\gamma_q(t) &=  x_0 + \sum_{j = 1}^2  q  v_j \exp(i [{\rm Imag}(\lambda_j) + \Delta \omega  {\rm sign}({\rm Imag}(\lambda_j))] t)   \nonumber \\
 &=x_0 + 2q \Big[ {\rm Real}(v_1) \cos( ({\rm Imag}(\lambda_1) + \Delta \omega)t)  - {\rm Imag}(v_1) \sin( ({\rm Imag}(\lambda_1) + \Delta \omega) t)  \Big], \nonumber \\
 \alpha(q,t) &=   \sum_{j = 1}^2 q  \Big[ - {\rm Real}(\lambda_j) + i {\rm sign}({\rm Imag}(\lambda_j)) \Delta \omega \Big] v_j \exp(i ({\rm Imag}(\lambda_j) + {\rm sign}({\rm Imag}(\lambda_j))) t) \nonumber \\
&=  -2 q {\rm Real}(\lambda_1) \Big[ {\rm Real}(v_1) \cos( ({\rm Imag}(\lambda_1) + \Delta \omega) t)    - {\rm Imag}(v_1) \sin( ({\rm Imag}(\lambda_1) + \Delta \omega) t)  \Big]   \nonumber \\
&\qquad - 2 \Delta \omega q \Big[ {\rm Imag}(v_1) \cos( ({\rm Imag}(\lambda_1) + \Delta \omega) t)  + {\rm Real}(v_1) \sin( ({\rm Imag}(\lambda_1) + \Delta \omega) t)   \Big]    .
 \end{align}
 The application of the input $\alpha(q,t)$ mandated by Equation \eqref{inporb} yields a periodic orbit $x^\gamma_q(t)$ with period $T = 2\pi/({\rm Imag}(\lambda_1) + \Delta \omega)$.  This orbit can be written in the form \eqref{yorb} and subsequently analyzed according to a greatly simplified version of the phase-amplitude reduction \eqref{forcedadaptive}.  To illustrate this, consider that the Jacobian associated with the periodic orbit \eqref{yorb} can be written as
 \begin{align} \label{jaceq}
 J(t) &=  \begin{bmatrix} J_0 & \frac{\partial \alpha}{\partial s} \\ 0&0\end{bmatrix} + \epsilon P(t) + O(\epsilon^2),
 \end{align}
 where $J_0$ is the Jacobian of $F$ evaluated at $x = x_0$ and  $u = 0$,    the zeros in Equation \eqref{jaceq} appropriately sized vectors,  $\frac{\partial \alpha}{\partial s}$ is evaluated at $s=t$, and $P(t)$ is comprised of second order partial derivatives of $F$.  With \eqref{jaceq} in mind, considering the periodic orbit of the form \eqref{yorb}, letting $\Delta y = y - y_q^\gamma(t)$, nearby solutions evolve to leading order according to
 \begin{equation} \label{liny}
 \Delta \dot{y} = \begin{bmatrix} J_0 & \frac{\partial \alpha}{\partial s} \\ 0&0\end{bmatrix} \Delta y + \epsilon P(t) \Delta y.
 \end{equation}
 Solutions of \eqref{liny} can be used to obtain Floquet eigenfunctions and Floquet eigenvalues associated with the periodic orbit $y^\gamma_q$.  Note that Equation \eqref{liny} is linear time varying system with period $T = 2\pi/({\rm Imag}(\lambda_1) + \Delta \omega)$.  As such, the monodromy matrix $\Phi$, i.e.,~that yields the relationship $\Delta y (T) = \Phi y(0)$ can be approximated by first noticing that \eqref{liny} is linear time invariant to leading order $\epsilon$ so that 
 \begin{align} \label{fundmtx}
 \Phi &=  \exp \bigg( \begin{bmatrix} J_0 & \frac{\partial \alpha}{\partial s} \\ 0&0\end{bmatrix} T\bigg) + O(\epsilon) \nonumber \\
 &=  \begin{bmatrix}  \exp(J_0 T) & X \\ 0 & 1  \end{bmatrix}  + O(\epsilon) \nonumber \\
 &= \Phi_0 + O(\epsilon),
 \end{align}
 where $\exp(\cdot)$ denotes the matrix exponential, $X \in \mathbb{R}^N$, and $\Phi_0$ is comprised of the $O(1)$ terms of $\Phi$.  Above, the second line is obtained by exploiting the block triangular structure of the $O(1)$ terms.

 Eigenvalues and corresponding left and right eigenvectors of $\Phi$ (denoted by $\lambda_j^\Phi$, $w_j^\Phi$, and $v_j^\Phi$, respectively) determine the Floquet multipliers and Floquet eigenfunctions, respectively, associated with the periodic orbit from \eqref{yorb}.  Note here that the notation $(\lambda_j^\Phi,v_j^\Phi)$ is used to denote the eigenvalue and eigenvector pair of the monondromy matrix $\Phi$; this is different from an eigenvalue and eigenvector pair of $J_0$ which is denoted by $(\lambda_j,v_j)$.   As discussed in Appendix \ref{apxb}, provided $\lambda_j^{\Phi_0}$ is simple, $O(\epsilon)$ perturbations to $\Phi_0$ will yield $O(\epsilon)$ shifts to the resulting eigenvalues and eigenvectors.  Thus, to leading order $\epsilon$, Floquet multipliers and Floquet eigenfunctions can be obtained by considering the eigenvalues and eigenvectors of $\Phi_0$.  For any eigenvalue, right eigenvector, and left eigenvector triple $(\lambda_j,v_j,w_j)$ of $J_0$
 \begin{equation} \label{initvec}
v_j^{\Phi_0}  = \begin{bmatrix} v_j  \\ 0 \end{bmatrix},
 \end{equation}
 is a right eigenvector of $\Phi_0$ with eigenvalue $\lambda_j^{\Phi_0} = \exp(\lambda_j T)$.   One can also verify that 
 \begin{equation} \label{initleft}
w_j^{\Phi_0} = \begin{bmatrix} w_j  \\  X^T w_j/(\exp(\lambda_j T) - 1)   \end{bmatrix},
 \end{equation}
  is an eigenvector of $\Phi_0^T$ with eigenvalue $\lambda_j^{\Phi_0} = \exp(\lambda_j T)$, and hence, is the corresponding left eigenvector of $\Phi_0$.  With this in mind, the corresponding Floquet exponent associated with the periodic orbit from \eqref{yorb} is
 \begin{align} \label{floqval}
 \kappa_j &= \frac{\log(\exp(\lambda_j T))}{T} + O(\epsilon) \nonumber \\
 &= \frac{\log | \exp(\lambda_j T) |}{T} +  i \frac{{\rm arg}(\exp(\lambda_j T))}{T} + O(\epsilon)  \nonumber \\
 &= {\rm Real}(\lambda_j) +  i \frac{{\rm arg}( {\rm exp} (i {\rm Imag}(\lambda_j) T  ))}{T} + O(\epsilon) \nonumber \\
 &= {\rm Real}(\lambda_j) + i {\rm Imag}(\lambda_j) - \frac{2 \pi i m }{T} + O(\epsilon). 
 \end{align}
 Above, note that $\exp(a + b i) = \exp(a + b i + 2 k \pi i)$ for any value of $k$.  As such, \eqref{floqval} mandates that the imaginary component of $\log(\exp(\lambda_j T))$ is always in the interval $(-\pi,\pi]$.  As such the value of $m$ must be chosen appropriately in the final line for each Floquet exponent.  Also notice that because of the constraints on $\Delta \omega$ given below Equation \eqref{levelsets}, $m=1$ when considering $\kappa_1$; this point will become important momentarily.   Associated Floquet eigenfunctions can be obtained by finding periodic solutions of Equation \eqref{adjfloq}.  Toward this end, considering the Jacobian from Equation \eqref{jaceq} to leading order $\epsilon$ Equation \eqref{adjfloq} becomes
\begin{equation} \label{preg}
 \dot{g}_j =   \bigg( \begin{bmatrix} J_0 & \frac{\partial \alpha}{\partial s} \\ 0&0\end{bmatrix}  - \kappa_j {\rm Id} \bigg) g_j + O(\epsilon).
\end{equation}
To leading order $\epsilon$, initial conditions that produce periodic solutions to \eqref{preg} are given by the right eigenvectors of $\Phi$ which were given in \eqref{initvec}.  With this in mind, notice that for any $g_j \propto \begin{bmatrix} v_j^T & 0\end{bmatrix}^T$, equation \eqref{preg} simplifies to
 \begin{align} \label{gic} 
 \dot{g}_j &= g_j  (\lambda_j - \kappa_j) + O(\epsilon) \nonumber \\
&= g_j   ({\rm Imag}(\lambda_j) - {\rm Imag}(\kappa_j))i  + O(\epsilon) \nonumber \\
 &=g_j   2 \pi i m /T + O(\epsilon) \nonumber \\
 &=g_j   ({\rm Imag}(\lambda_1) + \Delta \omega) i m  + O(\epsilon),
 \end{align}
where $m$ is an integer.  Above, the third line is obtained by substituting the final line of Equation \eqref{floqval} and the fourth line is obtained by noting that $T = 2 \pi/ ({\rm Imag}(\lambda_1) + \Delta \omega)$.  Considering Equation \eqref{gic}, to leading order $\epsilon$, if $g_j$ is proportional to $\begin{bmatrix} v_j^T & 0\end{bmatrix}^T$ at time $t = 0$, it remains proportional on timescales of order $1/\epsilon$.  As such, integrating Equation \eqref{gic} over one period ultimately yields the solution
 \begin{equation} \label{geq}
 g_i(t) = \begin{bmatrix} v_j \\ 0 \end{bmatrix} \exp( ( {\rm Imag}(\lambda_1) + \Delta \omega) i m t) + O(\epsilon).
 \end{equation}
The gradient of the Floquet coordinates can be obtained by finding solutions of Equation \eqref{adjiso}.  To leading order $\epsilon$, this equation is given by
\begin{equation} \label{irceq}
\dot{I}_j = -  \bigg( \begin{bmatrix} J_0^T & 0 \\ \frac{\partial \alpha}{\partial s}^T &0\end{bmatrix}  - \kappa_j {\rm Id} \bigg) I_j + O(\epsilon).
\end{equation}
 Similar to how solutions $g_j$ were obtained, using the left eigenvector obtained from Equation \eqref{initleft}, as an initial condition one can show that Equation \eqref{irceq} along solutions can be written as
 \begin{equation} \label{ieps}
 \dot{I}_j = \begin{bmatrix} -I_{j,1}   ({\rm Imag}(\lambda_1) + \Delta \omega) i m  \\ -( \frac{\partial \alpha}{\partial s}^T I_{j,1} - \kappa_j I_{j,2} ) \end{bmatrix} + O(\epsilon),
 \end{equation}
 where $I_{j,1}(t) \in \mathbb{C}^N$ are the first $N$ elements of $I_j(t)$ and $I_{j,2}(t) \in \mathbb{C}$ is the last element.   As such, solutions of \eqref{irceq} have the following form:
  \begin{equation} \label{ieq}
 I_j(t) = \begin{bmatrix}  w_j  \exp(-({\rm Imag}(\lambda_1) + \Delta \omega) i m t)  \\ \rho_j(t)  \end{bmatrix}  + O(\epsilon),
 \end{equation}
 where $\rho_j(t) \in \mathbb{C}$ is periodic.   Further simplification for  $\rho_1(t)$ from Equation \eqref{ieq} is possible by considering the constraint \eqref{normeq} from Appendix \ref{apx0}.  With a change of variables $\theta = \omega t$, this constraint becomes
 \begin{equation} \label{mainconst}
 I_1^T(t) \frac{\partial y^\gamma}{\partial t} = 0.
 \end{equation}
The term  $\partial y^\gamma/\partial \theta$ can be obtained from direct differentiation of \eqref{inporb}; substituting into \eqref{mainconst} one finds
 \begin{align}
 0 &= I_1^T(t) \frac{\partial y^\gamma}{\partial t} \nonumber \\
 &= \begin{bmatrix} w_1^T  \exp(- ({\rm Imag} (\lambda_1) + \Delta \omega) i  t ) & \rho(t) \end{bmatrix}   \nonumber \\
 & \quad \times  \begin{bmatrix}  \sum_{j = 1}^2 q i  ({\rm Imag}(\lambda_j) + {\rm sign}({\rm Imag}(\lambda_j))\Delta \omega) v_j \exp( ({\rm Imag}(\lambda_j)+ {\rm sign}({\rm Imag}(\lambda_j))\Delta \omega) i t) \\  1 \end{bmatrix} \nonumber \\
 &= q i ({\rm Imag}(\lambda_1) + \Delta \omega) + \rho_1(t).
 \end{align}
 In the second line above,  $m = 1$ for $I_1$ as mentioned directly below Equation \eqref{floqval} yielding the simplification in the final line.    Ultimately, one finds $\rho_1(t) = -q i ({\rm Imag}(\lambda_1) + \Delta \omega)$ so that
  \begin{equation} \label{ieq}
 I_1(t) = \begin{bmatrix}  w_1  \exp(-({\rm Imag}(\lambda_1) + \Delta \omega) i t)  \\ -q i ({\rm Imag}(\lambda_1) + \Delta \omega)  \end{bmatrix}  + O(\epsilon).
 \end{equation}

 \subsection{Adaptive Phase Amplitude Reduction For an Oscillatory Mode Near a Fixed Point} \label{adaptsing}
Considering the derivations for the input $\alpha(q,t)$ and the gradient of the Floquet coordinates $I_j(t)$  from Section \ref{singmode}, the adaptive reduction from Equation \eqref{forcedadaptive} admits substantial simplifications.   Note that in Equations \eqref{inporb}, \eqref{geq}, and \eqref{ieq}, $\theta = \frac{2 \pi}{T} t = \omega t$ along trajectories.  As such, one can use the substitution $t = \theta/\omega$ when considering phase coordinates.  For the purposes of this subsection, it will be assumed that the adaptive parameter, $q$, is an $O(\epsilon)$ term.  The dynamics of $\psi_1$ as mandated by Equation \eqref{forcedadaptive} are
\begin{align} \label{psi1int}
\dot{\psi_1} &= \kappa_1(q) \psi_1 + I_{1,1}^T(\theta,q)U_e + I_{1,2}(\theta,q) f_\theta + E_1(\theta,q) \dot{q} \nonumber \\
&= \kappa_1(q) \psi_1 +  \exp(- i \theta) w_1^T U_e - q i \omega f_\theta + E_1(\theta,q) \dot{q},
\end{align} 
 where $I_{1,1}(\theta,1)$ and $I_{1,2}$ were defined below Equation \eqref{ieps} as the first and last components of $I_1(\theta,q)$.  In the second line, the relation from Equation \eqref{ieq} is used, the terms are written as functions of $\theta$ using the coordinate transformation $\theta = \omega t$ where $\omega = ({\rm Imag}(\lambda_1) + \Delta \omega)$.  Further, as discussed in Appendix \ref{apx0}, $E_1(\theta,q)$ can be computed according to
 \begin{align}
 E_1(\theta,q) &= -I_1^T(\theta,q) \frac{\partial x^\gamma}{ \partial q} \nonumber \\
 &= - \begin{bmatrix}  w_1 \exp(-i \theta) \\ -q i {\rm Imag}(\lambda_1) \end{bmatrix}^T  \begin{bmatrix} \sum_{j = 1}^2  v_j \exp((-1)^{j+1} i \theta) \\ 0 \end{bmatrix} \nonumber \\
 &= -1,
 \end{align}
where the second line is obtained by substituting $\eqref{ieq} $ for $I_1(\theta,q)$ and using \eqref{inporb} to obtain $\partial x^\gamma/\partial q$ as defined in Appendix \ref{apx0}.  In Equation \eqref{psi1int}, $\dot{q}$ and $f_\theta$ can be chosen arbitrarily; specifically taking these terms to be 
\begin{align}
\dot{q} &= {\rm Real}(\exp(- i \theta) w_1^T) U_e, \nonumber \\
f_\theta &= \frac{ {\rm Imag}(\exp(- i \theta) w_1^T) U_e}{q \omega},
\end{align}
yields
\begin{equation}
\dot{\psi}_1 = \kappa_j(q) \psi_1.
\end{equation}
Recalling from Equation \eqref{floqval} that ${\rm Real}(\kappa_j(q)) = {\rm Real}(\lambda_j)   + O(\epsilon)$ when $q$ is an order $\epsilon$ term, ${\rm Real}(\kappa_j(q))<0$ for $q$ small enough so that $\lim_{t \rightarrow \infty} \psi_1 = 0$ and the dynamics of both $\psi_1$ and $\psi_2 = \psi_1^*$ can be ignored.  Taken together, the adaptive reduction from Equation \eqref{forcedadaptive} becomes
\begin{align} \label{smallqadapt}
\dot{\theta} &= \omega(q) + \frac{   {\rm Imag}(\exp(- i \theta) w_1^T) U_e}{q}, \nonumber \\
\dot{q} &=  {\rm Real}(\exp(- i \theta) w_1^T) U_e, \nonumber \\
\dot{\psi}_j &= \kappa_j(q) \psi_j + I^T_{j,1}(\theta,q)U_e + I_{j,2}(\theta,q) f_\theta + E_j(\theta,q) \dot{q}, \nonumber \\
j & = 3, \dots, \beta.
\end{align}
Once again, Floquet coordinates $\psi_j$ are ignored if the corresponding value of ${\rm Real}(\kappa_j(q))$  is negative and large in magnitude so that transients decay rapidly.

It is worthwhile to consider the dynamics of Equation \eqref{smallqadapt} in the unperturbed setting, i.e.,~when $u = 0$.  In this case, $U_e = -\alpha(q,\theta)$, which can be obtained from Equation \eqref{inporb} after substituting $\theta = ({\rm Imag}(\lambda_1) + \Delta \omega) t$.   Directly multiplying these terms together and simplifying (with help from the fact that $w_1^T v_1 = 1$ and $w_1^T v_1^* = 0$ so that ${\rm Real}(w_1^T) {\rm Real}(v_1) = 0.5$, ${\rm Imag}(w_1^T) {\rm Imag}(v_1) = -0.5$, ${\rm Real}(w_1^T) {\rm Imag}(v_1) = 0$, and ${\rm Imag}(w_1^T) {\rm Real}(v_1) = 0$) one finds that when $u = 0$, the $\theta$ and $q$ dynamics of \eqref{smallqadapt} simplify to
\begin{align} \label{actionangle}
\dot{\theta} &= {\rm Imag}(\lambda_1), \nonumber \\
\dot{q} &= q {\rm Real}(\lambda_1).
\end{align}
Note that the representation from Equation \eqref{actionangle} gives dynamics that are identical to those of an action-angle coordinate system (see for instance \cite{maur13}) that are valid for a linear system with a complex eigenvalue.  Note that the results from Sections \ref{singmode} and \ref{adaptsing} are valid when $q = O(\epsilon)$, i.e.,~when the state is close enough to the underlying fixed point so that the dynamics can be represented according to a local linearization.  As illustrated in the following sections, however, it is possible to extend this strategy to consider nonlinear oscillations far beyond a close neighborhood of the stable fixed point.

\subsection{Extension For Nonlinear Oscillations Far Beyond the Stable Fixed Point} \label{nlinmode}
Sections \ref{singmode} and \ref{adaptsing} consider the construction of appropriate trajectories for the adaptive reduction for locations close to the fixed point.  This strategy can be readily extended to consider trajectories far beyond the fixed point.  Do do so, first consider a general periodic orbit $y_{q_0}^\gamma$, of Equation \eqref{nonaut} that results when $u = 0$ with $q_0 \in \mathbb{R}$ being an $O(\epsilon)$ term.  Recall that this is a forced periodic orbit with $\alpha(q_0,t)$ of the form given in Equation \eqref{inporb}.   As in Section \eqref{singmode}, we let $g_1(t)$ and $g_2(t)$ be Floquet eigenmodes associated with the Floquet exponents $\kappa_1$ and $\kappa_2$ as given in Equations \eqref{geq} and \eqref{floqval}, respectively.  One can define an adjacent periodic orbit
\begin{align}  \label{xdelt}
y^\gamma_{q_0 + \Delta q}(t) &= y^\gamma_{q_0}(t) + \sum_{j = 1}^2 \Delta {q} g_j(t) \nonumber \\
&=y^\gamma_{q_0}(t) + 2 \Delta q {\rm Real}(g_1(t)),
\end{align}
where $\Delta q$ is also an $O(\epsilon)$ term and the simplification in the second line results from the fact that $g_1(t) = g_2^*(t)$.  Next, we seek external forcing $\alpha(q_0 +\Delta q,t) = \alpha(q_0,t) + \Delta \alpha(q_0,t)$ that admits the periodic orbit mandated by Equation \eqref{xdelt}.  To proceed, taking the time derivative of \eqref{xdelt} yields
\begin{equation} \label{dotx1}
\dot{y}^\gamma_{q_0 + \Delta q}(t) = \frac{d}{dt} \begin{bmatrix}  x^\gamma_{q_0 + \Delta q}(t) \\ {\rm mod}(t,T) \end{bmatrix} = \dot{y}^\gamma_{q_0}(t) +  \sum_{j = 1}^2 \Delta q \dot{g}_j(t).
\end{equation}
Likewise, along this new periodic orbit, by considering the $x$ dynamics of the underlying model \eqref{nonaut} one also finds
\begin{align} \label{dotx2}
\dot{x}^\gamma_{q_0 + \Delta q}(t)   &= F({x}^\gamma_{q_0 + \Delta q}(t) ,0) + \alpha(q_0,t) + \Delta \alpha(q_0,t), \nonumber \\
&= F({x}^\gamma_{q_0}(t) ,0)  + \frac{\partial F}{\partial x} \bigg( \sum_{j = 1}^2 \Delta q {g}_{j,1}(t) \bigg) + \alpha(q_0,t) + \Delta \alpha(q_0,t) + O(\epsilon^2).
\end{align}
In the above equation, the partial derivatives are evaluated at $x^\gamma_{q_0}(t)$ and $g_{j,1}(t) \in \mathbb{C}^N$ corresponds to the first $N$ elements of $g_j$.  Combining Equations \eqref{dotx2} and the first $N$ terms of  \eqref{dotx1}, noting that $\dot{x}^\gamma_{q_0}(t)  = F({x}^\gamma_{q_0}(t) ,0)  + \alpha(q_0,t)$ one finds
\begin{align} \label{deltaalpha}
\Delta \alpha(q_0,t) &= \sum_{j = 1}^2 \Delta q \dot{g}_{j,1}(t) - \frac{\partial F}{\partial x} \bigg( \sum_{j = 1}^2 \Delta q {g}_{j,1}(t) \bigg) + O(\epsilon^2) \nonumber \\
&= 2 \Delta q {\rm Real}(\dot{g}_{1,1}(t)) - 2 \Delta q  \frac{\partial F}{\partial x} {\rm Real}(g_{1,1}(t)) + O(\epsilon^2),
\end{align}
where the second line can be obtained by noting that $g_1(t) = g_1^*(t)$.  As such, choosing
\begin{equation} \label{alphadelt}
\alpha(q_0 + \Delta q,t) = \alpha(q_0,t) + 2 \Delta q {\rm Real}(\dot{g}_{1,1}(t)) - 2 \Delta q  \frac{\partial F}{\partial x} {\rm Real}(g_{1,1}(t)),
\end{equation}
yields the periodic orbit mandated by \eqref{xdelt} to leading order $\epsilon$.  

One can verify that to leading order $\epsilon$, Equations \eqref{deltaalpha} and \eqref{alphadelt} yield the same results as given in Equation \eqref{inporb} (recalling that $q_0$ and $\Delta q$ are both $O(\epsilon)$ terms).   Nonetheless, these definitions provide a strategy for extending the proposed method beyond the linear regime.   To this end, suppose that an orbit $x^\gamma_{q_0 + k \Delta q}$ has already been obtained with associated input $\alpha(q_0 + k \Delta q,t)$ where $k \in \mathbb{N}$.  Both of the terms $I_j(t,q_0 + k \Delta q)$ and $g_j(t,q_0 + k \Delta q)$ can be obtained from Equations \eqref{adjiso} and \eqref{adjfloq}, respectively.  Here, we note that $q_0 + k \Delta q$ is no longer constrained to be an $O(\epsilon)$ term so that $g_j$ and $I_j$ are now functions of both $t$ and $q$.  Each $g_j(q,t)$ must be normalized appropriately so that $g(t,q)$ is continuous; recall that the eigenvectors $v_1$ and $v_2$ defined above Equation \eqref{linmode1} are scaled so that $||v_1||_2 = ||v_2||_2=1$, $v_1 = v_2^*$, and ${\rm arg}(e_j^T v_1) = -\pi$ where $e_j$ is the $j^{\rm th}$ element of an appropriately sized standard unit basis and $j$ can be chosen arbitrarily. In order to match the scaling on  $g_1(t,q)$ and $g_2(t,q)$ resulting from Equation \eqref{geq}, we require that $||g_1(0,q)||_2 = ||g_2(0,q)||_2 = 1$, $g_1(0,q) = g_2^*(0,q)$ and ${\rm arg}(e_j^T g_1(0,q)) = -\pi$ for all $q$.  Subsequently, a new periodic orbit and associated external input can be defined according to
\begin{align} \label{nlinupdate}
x^\gamma_{q_0 + (k+1) \Delta q}(t) &= x^\gamma_{q_0 + k \Delta q}(t) + \Delta q \sum_{j = 1}^2 g_i(t,q_0 + k \Delta q), \nonumber \\
\alpha(q_0 + (k+1)\Delta q,t) &= \alpha(q_0+k\Delta q,t) + 2 \Delta q {\rm Real}(\dot{g}_{1,1}(t,q_0 + k\Delta q)) \nonumber \\
& \quad -2 \Delta q \frac{\partial F}{\partial x} {\rm Real}(g_{1,1}(t,q_0 + k \Delta q)).
\end{align}
Note that the periodic orbit defined by \eqref{nlinupdate} is accurate to leading order in $\Delta q$.  For the purposes of practical implementation, in order to avoid compounding these errors at each iteration it is generally necessary to find the nearby, truly periodic solution which can be accomplished, for instance, using a Newton iteration.  This process can be repeated to define a family of periodic orbits that extend beyond the linear regime where a local linearization would be valid.  For this family of periodic orbits, considering the adaptive reduction from Equation \eqref{forcedadaptive}, the dynamics of the $\psi_1$ Floquet coordinate are 
\begin{equation} \label{psi1eq}
\dot{\psi}_1  =\kappa_1(q) \psi_1 + I_{1,1}^T(\theta,q) U_e + I_{1,2}(\theta,q) f_\theta + E_1(\theta,q) \dot{q}.
\end{equation}
As in Section \ref{adaptsing}, considering $E_1(\theta,q)$ as discussed in Appendix \ref{apx0}, $E_1(\theta,q)$ can be computed according to
 \begin{align}
 E_1(\theta,q) &= -I_1^T(\theta,q) \frac{\partial x^\gamma}{ \partial q} \nonumber \\
 &= -   I_1^T(\theta,q)   \bigg( \sum_{j = 1}^2 g_i(\theta , q_0 + k \Delta q) \bigg) \nonumber \\
 &= -1,
 \end{align}
where the second line is obtained by taking the partial of $\partial x^\gamma/\partial q$ from \eqref{nlinupdate} and interpreting each $g_i$ as a function of $\theta$ using the change of variables $\theta = \omega t$; the third line follows using the relation \eqref{gconst}.  Once again, recalling that $f_\theta$ and $\dot{p}$ from the adaptive reduction \eqref{forcedadaptive} can be chosen arbitrarily.  When taking
\begin{align} \label{inveq}
\begin{bmatrix}  \dot{q} \\ f_\theta  \end{bmatrix}  &= \begin{bmatrix} 1 & -{\rm Real}(I_{1,2}(\theta,q)) \\ 0 & -{\rm Imag}(I_{1,2}(\theta,q)) \end{bmatrix}^{-1} \begin{bmatrix}  {\rm Real}(I^T_{1,1}(\theta,q)) U_e \\  {\rm Imag}(I^T_{1,1}(\theta,q)) U_e \end{bmatrix},  \nonumber \\
\end{align}
Equation \eqref{psi1eq} becomes $\dot{\psi}_1 = \kappa_1(q) \psi_1$ which tends to zero in the limit as $t$ approaches infinity provided ${\rm Real}(\kappa_j(q)) < 0$ for all allowable $q$.  Using Equation \eqref{inveq} to determine $\dot{q}$ and $f_\theta$, the dynamics of both $\psi_1$ and $\psi_2=\psi_1^*$ can be ignored.  For the inverse in Equation \eqref{inveq} to exist, it is necessary that ${\rm Imag}(I_{1,2}(\theta,q)) \neq 0$; recall that this is guaranteed in the limit that $q$ is small as discussed in Section \ref{singmode}.  

Similar to the result from Section \ref{adaptsing}, taking $\dot{q}$ and $f_\theta$ as mandated by \eqref{inveq} the adaptive reduction from \eqref{forcedadaptive} becomes
\begin{align} \label{largeqadapt}
\dot{\theta} &= \omega(q) + \omega(q) f_\theta, \nonumber \\
\dot{q} &=  G_q(q,\theta,U_e), \nonumber \\
\dot{\psi}_j &= \kappa_j(q) \psi_j + I^T_{j,1}(\theta,q)U_e + I_{j,2}(\theta,q) f_\theta + E_j(\theta,q) \dot{q}, \nonumber \\
j & = 3, \dots, N,
\end{align}
where $G_q(q,\theta,U_e) = {\rm Real}(I^T_{1,1}(\theta,q)) U_e   -  {\rm Real}(I_{1,2}(\theta,q))  {\rm Imag}(I^T_{1,1}(\theta,q)) U_e /{\rm Imag}(I_{1,2}(\theta,q))$ and $f_\theta = -{\rm Imag}(I^T_{1,1}(\theta,q)) U_e / {\rm Imag}(I_{1,2}(\theta,q))$.   In the limit that $q$ is small, Equation \eqref{largeqadapt} reduces to Equation \eqref{smallqadapt}.  Nonetheless, the family of periodic orbits defined iteratively by \eqref{nlinupdate} can extended to regimes where simple linearization techniques are no longer valid.  

In the results to follow, the effective unforced natural frequency will also be considered, defined as the unforced natural frequency averaged over all $\theta$ that results when taking $u = 0$.  In this case, $U_e = -\alpha(q,\theta)$ so that 
\begin{align} \label{avfreq}
\bar{\omega}(q) &= \frac{1}{2\pi} \int_0^{2\pi} \big[ \omega(q) + \omega(q) f_\theta \big]d\theta \nonumber \\
&=  \omega(q)  +  \frac{1}{2\pi} \int_0^{2\pi} \bigg[ \omega(q)  \frac{ {\rm Imag}(I^T_{1,1}(\theta,q)) \alpha(q,\theta) }{ {\rm Imag}(I_{1,2}(\theta,q)) }  \bigg]  d\theta, \nonumber \\
\end{align}
where $f_\theta$ is taken as mandated by Equation \eqref{inveq}.  Considering the definition from \eqref{avfreq} and the result from Equation \eqref{actionangle} which is valid when $q = O(\epsilon)$  (i.e.,~when the state is close to the underlying fixed point), $\bar{\omega}(q) = {\rm Imag}(\lambda_1)$ to leading order $\epsilon$ when $q$ is an order $\epsilon$ term.  

\subsection{Considering External Forcing with a Nonstatic Period} \label{nonstatper}

The family of periodic orbits constructed iteratively by Equation \eqref{nlinupdate} all have the same period.  In some cases, it may be desirable to change the natural frequency of a given periodic orbit.  To this end, consider a $T$-periodic solution $x^\gamma_{q}(t)$  that provides a solution to Equation \eqref{fulleq} when applying the  $T$-periodic input  $\alpha(q,t)$ and taking $u = 0$.  Recalling that $\omega = 2\pi/T$, the periodic orbit $x^\gamma_q(    \frac{\omega + \Delta \omega}{\omega} t)$ can be obtained instead by making an appropriate change to $\alpha(q,t)$.  To this end, substituting the desired periodic orbit into Equation \eqref{fulleq} yields
\begin{align} \label{findalpha}
\dot{x}^\gamma_q \bigg(    \frac{\omega + \Delta \omega}{\omega} t \bigg)  \bigg( 1 +   \frac{\Delta \omega}{\omega} \bigg) = F\bigg( x^\gamma_q \bigg(    \frac{\omega + \Delta \omega}{\omega} t \bigg),0 \bigg) + \hat{\alpha}(q,t),
\end{align}
where $\hat{\alpha}$ is the input which provides a solution to Equation \eqref{findalpha}.  Noting that $\dot{x}^\gamma_q =F(x_q^\gamma(t),u) + \alpha(q,t)$, taking
\begin{equation} \label{alphaeq}
\hat{\alpha}(q,t) = \alpha\bigg(q, \bigg( \frac{\omega + \Delta \omega}{\omega} \bigg) t \bigg) + \frac{\Delta \omega}{\omega} \dot{x}^\gamma_q \bigg(    \frac{\omega + \Delta \omega}{\omega} t \bigg),
\end{equation}
provides a solution for Equation \eqref{findalpha} and yields the desired periodic orbit.  Considering the autonomous system of equations from Equation \eqref{nonaut}, recalling that $s\in[0,T)$ one finds
\begin{align} \label{xshift}
\dot{x} &= F(x,u) + \hat{\alpha}(q,s), \nonumber \\
\dot{s} &= \frac{\omega + \Delta \omega}{\omega},
\end{align}
admits the appropriately time shifted periodic orbit $x^\gamma_q(  \frac{\omega + \Delta \omega}{\omega} t)$ when $u = 0$ with period $\hat{T}=2\pi/(\omega + \Delta \omega)$ and natural frequency $\hat{\omega} = \omega + \Delta \omega$.  This shift in the externally applied forcing can be implemented in conjunction with the  iteratively defined periodic orbits given in Equation \eqref{nlinupdate}.  Noting that $\partial x^\gamma/\partial q$ remains unchanged when the period of a given orbit is shifted, when considering nonstatic periods the general structure of the adaptive reduction from \eqref{largeqadapt} remains unchanged with the main difference being that $\omega$ becomes a function of $q$.  

Changing the period of oscillation according to the strategy above has the general effect of changing the Floquet exponents associated with the underlying periodic orbit.  For instance, when $q = O(\epsilon)$, consider a shift in the natural period from $T$ to $T + \Delta T$  that results from the modification to the external forcing  mandated by Equation \eqref{alphaeq}.  Despite the shift in period, because the underlying orbit remains unchanged, comparing to Equation \eqref{fundmtxshift}, the monodromy matrix associated with the periodic orbit of \eqref{xshift} can be obtained according to
 \begin{align} \label{fundmtxshift}
 \Phi &=  \exp \bigg( \begin{bmatrix} J_0 & \frac{\partial \alpha}{\partial s} \\ 0&0\end{bmatrix} T\bigg) + O(\epsilon) \nonumber \\
 &=  \begin{bmatrix}  \exp(J_0 T) & \hat{X} \\ 0 & 1  \end{bmatrix}  + O(\epsilon) \nonumber \\
 &= \hat{\Phi}_0 + O(\epsilon),
 \end{align}
where the term $\hat{X}$ above differs from $X$ in \eqref{fundmtx} because of the difference in $\partial \alpha/\partial s$.  Nonetheless, the right eigenvector from \eqref{initvec} is also a right eigenvector of $\hat{\Phi}$ with eigenvalue $\exp(\lambda_j (T + \Delta T))$.  Considering the relationship from \eqref{floqval}, the corresponding Floquet exponent is
\begin{align} \label{floqchange}
\hat{\kappa}_j &= \frac{\log(\exp(\lambda_j (T+\Delta )))}{T + \Delta T} + O(\epsilon) \nonumber \\
&=  \frac{\log | \exp(\lambda_j (T + \Delta T)) |}{T + \Delta T} +  i \frac{{\rm arg}(\exp(\lambda_j (T + \Delta T)))}{T+\Delta T} + O(\epsilon) \nonumber \\
&= {\rm Real} (\lambda_j) +  i \frac{{\rm arg}(\exp(i {\rm Imag} (\lambda_j) (T + \Delta T)))}{T+\Delta T}+ O(\epsilon) \nonumber \\
&= {\rm Real}(\lambda_j) + i \frac{{\rm Imag}(\lambda_j)(T + \Delta T) - 2 \pi m}{T+\Delta T} + O(\epsilon) \nonumber \\
&= \kappa_j + \frac{2 \pi m i \Delta T}{T^2} + O(\Delta T^2) + O(\epsilon).
\end{align}
Above, the final line is obtained with a Taylor expansion centered at $\Delta T = 0$.  Note that the relation \eqref{floqchange} is only valid in the limit that $q$ is an $O(\epsilon)$ term.  Comparing to the value of $\kappa_j$ from Equation \eqref{floqval},  changing the period will only change the imaginary component of the resulting Floquet exponent to leading order $\epsilon$. 

\subsection{List of Steps to Implement the Proposed Model Order Reduction Approach} \label{steplist}
A list of steps required to implement the proposed reduced order modeling strategy detailed in Sections \ref{overviewsec}-\ref{nonstatper} are summarized below:
\begin{enumerate}[1)]
\item Identify an appropriate fixed point $x_{\rm ss}$ of the model \eqref{maineq}, i.e.,~for which $F(x_{\rm ss},0) = 0$.  Identify a complex-conjugate pair of simple eigenvalues and eigenvectors $(\lambda_1,v_1)$ and $(\lambda_2,v_2)$ associated with a single oscillatory linear mode.  Order the eigenvalues so that ${\rm Imag}(\lambda_1)>0$.

\item Normalize $v_1$ and $v_2$ appropriately so that $||v_1||_2=||v_2||_2=1$, $v_1=v^*$ and ${\rm arg}(e_j^T v_1) = -\pi$.  Note that alternative normalizations could be used as long as they are used consistently in the implementation of this strategy.

\item Let $T = 2 \pi / ({\rm Imag}(\lambda_1) + \Delta \omega)$ where $\Delta \omega$ be a nonzero constant defined below Equation \eqref{levelsets}.  Also let $q_0$ be a small positive constant.  Obtain the $T$-periodic orbit $x_{q_0}^\gamma(t)$ that results when applying $\alpha(q_0,t)$, both found according to Equation \eqref{inporb}.  For this initial orbit, $q_0$ must be chosen small enough so the underlying dynamical system is well approximated by linearization.

\item  Considering the autonomous system represented according to \eqref{nonaut}, compute the terms $g_1(\theta,q_0)$ and $I_1(\theta,q_0)$ for the periodic orbit $y_{q_0}^\gamma(t)$ with terms that are computed in the previous step.  This can be accomplished by first finding the eigenvalue $\lambda_1^{\Phi_0} \approx \exp(\lambda_j T)$ associated with the monodromy matrix with associated Floquet exponent $\kappa_1(q_0)$ defined in Equation \eqref{floqval}.  The term $g_1(\theta,q_0)$ can be found by finding the periodic solution to Equation \eqref{adjfloq} and normalizing so that $||g_1(0,q_0)||_2 = 1$, and ${\rm arg}(e_j^T g_1(0,q_0)) = -\pi$ for all $q$ (which matches the normalization given in step 2 above).  Subsequently, one can compute $I_1(\theta,q_0)$ by finding the periodic solution to Equation \eqref{adjiso} using the normalization mandated by Equation \eqref{gconst}.

\item For a small positive constant $\Delta q$, define an adjacent periodic orbit  $x_{q_0 + \Delta q}^\gamma(t)$ with associated input $\alpha(q_0+\Delta q)$, both computed according to Equation \eqref{nlinupdate}.

\item Use the periodic orbit obtained in step 5, to define the periodic orbit $y_{q_0+\Delta q}^\gamma(t)$ for the autonomous system \eqref{nonaut}.  Compute $g_1(\theta,q_0+\Delta q)$ and $I_1(\theta,q_0+\Delta q)$ for this periodic orbit with appropriate scaling.  

\item Continue to iteratively define adjacent periodic orbits and associated inputs according to Equation \eqref{nlinupdate}.  For each newly identified orbit, compute the terms $g_1(\theta,q)$ and $I_1(\theta,q)$ for the associated periodic orbits. Note that the Floquet exponent $\kappa_1(q)$ will generally change slowly as the periodic orbits become farther from the fixed point.  The resulting information is used to define the dynamics governing $\theta$ and $q$ from Equation \eqref{largeqadapt}.
\end{enumerate}

We also emphasize a few general notes about the implementation of the proposed strategy below:
\begin{itemize}
\item The dynamics of a Floquet coordinate $\psi_j$ for $j \geq 3$ can generally be ignored if the associated Floquet exponent $\kappa_j(q)$ is negative and large in magnitude for all relevant $q$.  If a given isostable coordinate $\psi_j$ cannot be ignored, the terms $I_j(\theta,q)$ and $g_j(\theta,q)$ can be computed immediately after $I_1(\theta,q)$ and $g_1(\theta,q)$ for each periodic orbit.  Note that $\kappa_j(q) = \log(\lambda_j^\Phi)/T$ where $\lambda_j^\Phi$ an appropriate eigenvalue of the monodromy matrix.

\item As mentioned in the main text, when determining $\alpha(q + \Delta q,t)$ that yields the periodic orbit $x_{q+\Delta q}^\gamma(t)$, Equation \eqref{nlinupdate} is only valid up to order $\Delta q^2$.  In order to prevent these errors from accumulating over multiple iterations (and hence yielding non-periodic solutions) it is generally necessary to view each $x^\gamma_{q + \Delta q}(t)$ computed according to \eqref{nlinupdate} as a close guess and numerically identify a the periodic solution using a Newton iteration.

\item  As the value of $q$ increases, it is possible for $\lambda_1^\Phi$ to become a repeated eigenvalue if ${\rm Imag}(\lambda_1)T$ becomes a multiple of $2\pi$.  If this occurs, it is not guaranteed that the Floquet eigenfunction $g_1(\theta,q)$ will remain continuous with respect to $q$, thereby precluding the use of the adaptive phase amplitude reduction.  Slight modifications to the period using the strategy discussed in Section \ref{nonstatper} can be useful to prevent this issue.  Once $\alpha(q,t)$ and $x_{q+\Delta q}^\gamma(t)$ are computed for a given value of $q$, it is possible to adjust the period of oscillation by adjusting the external input according to Equation \eqref{alphaeq}. 
\end{itemize}

\subsection{Considering Multiple Adaptive Parameters for Multiple Nonlinear Oscillation Modes} \label{multimodes}
The analysis provided in Sections \ref{singmode}-\ref{nlinmode} yields a single adaptive parameter associated with oscillations of a single nonlinear mode.  Multiple oscillatory modes can also be considered with this formulation, but it is necessary for the magnitudes of associated Floquet coordinates  $\psi_3, \dots, \psi_N$ to remain small, i.e.,~the amplitudes of the other oscillatory modes must be small.  It is relatively straightforward to consider the influence of multiple nonlinear modes by considering additional adaptive parameters that ultimately yield an $n$-dimensional family of periodic orbits with $n>1$.  

To this end, suppose that a 1-parameter family of periodic orbits $x^\gamma_{q_1}(t)$ and associated input $\alpha(q_1,t)$ valid for $q_1 \in [0,q_{1,{\rm max}}]$ that has been obtained iteratively according to Equation \eqref{nlinupdate}.  In order to accommodate an additional nonlinear mode, one can use this initial set of periodic orbits to define a 3-parameter family or orbits $x^\gamma_{q_1,q_2,q_3}$ with an associated 3-parameter family of inputs $\alpha(q_1,q_2,q_3,t)$.  Similar to the strategy discussed in Section \ref{nlinmode}, let $g_3(q_1,q_2,q_3,t)$ and  $g_4(q_1,q_2,q_3,t)$ be Floquet eigenfunctions that correspond to complex-conjugate Floquet coordinates $\psi_3 = \psi_4^*$, i.e.,~corresponding to the second oscillatory mode (note that $g_1(q_1,q_2,q_3,t)$ and $g_2(q_1,q_2,q_3,t)$ are still the Floquet eigenfunctions associated with the first oscillatory mode).  Similar to the scaling on $g_1$ and $g_2$, we require $||g_3(q_1,q_2,q_3,0)||_2 = ||g_4(q_1,q_2,q_3,0)||_2 = 1$, $g_3(q_1,q_2,q_3,0) = g_4^*(q_1,q_2,q_3,0)$ and ${\rm arg}(e_j^T g_3(q_1,q_2,q_3,0)) = -\pi$ for all $q_1,q_2,q_3$ where $e_j$ is an appropriately sized element of the standard unit basis with $j$ chosen arbitrarily. 

 Starting by taking $x^\gamma_{q_1,0,0}(t)$ and $\alpha(q_1,0,0,t)$ to be identical to the periodic orbits and inputs defined obtained from Equation \eqref{nlinupdate} for a given value of $q_1$, a second dimension can be added as follows:
\begin{align}  \label{r1eq}
x^\gamma_{q_1,q_2+\Delta q_2,0}(t) &= x^\gamma_{q_1,q_2,0}(t) +   \Delta q_2 \sum_{j = 3}^4 g_j(q_1,q_2,0,t),
\end{align}
with corresponding input
\begin{align}
\alpha(q_1,q_2 + \Delta q_2,0,t) &= \alpha(q_1,q_2,0,t) +  \Delta q_2 \sum_{j = 3}^4  \dot{g}_j(q_1,q_2,0,t) -  \frac{\partial F}{\partial x}\bigg( \sum_{j = 3}^4  \Delta q_2  g_j(q_1,q_2,0,t) \bigg),
\end{align}
for $q_2 \in [q_{2,{\rm min}},q_{2,{\rm max}}]$.  Subsequently, a third dimension can be added according taking
\begin{equation} \label{r2eq}
x^\gamma_{q_1,q_2,q_3+\Delta q_3}(t) = x^\gamma_{q_1,q_2,q_3}(t) +   i \Delta q_3 \sum_{j = 3}^4 g_j(q_1,q_2,q_3,t),
\end{equation}
with corresponding input
\begin{align} \label{finalalpha}
\alpha(q_1,q_2,q_3 + \Delta q_3,t) &= \alpha(q_1,q_2,q_3,t) +  i \Delta q_3 \sum_{j = 3}^4  \dot{g}_j(q_1,q_2,q_3,t) -  \frac{\partial F}{\partial x}\bigg( \sum_{j = 3}^4  i \Delta q_3  g_j(q_1,q_2,q_3,t) \bigg),
\end{align}
for $q_3 \in [q_{3,{\rm min}},q_{3,{\rm max}}]$.  Note that similar to the periodic orbits and inputs defined iteratively according to Equation \eqref{nlinupdate}, the periodic orbits defined by Equations \eqref{r1eq} and \eqref{r2eq} are accurate to leading order in $\Delta q$.  For the purposes of practical implementation, in order to avoid compounding these errors at each iteration it is generally necessary to find the nearby, truly periodic solution which can be done, for instance, by using a Newton iteration.

Considering the additional dimensions for the family of periodic orbits defined by Equations \eqref{r1eq} and \eqref{r2eq}, letting $\vec{q} \equiv [q_1,q_2,q_3]^T$, the dynamics of the $\psi_1$ and $\psi_3$ Floquet coordinates are 
\begin{align} \label{psi1psi2}
\dot{\psi}_1 &= \kappa_1(\vec{q}\,) \psi_1 + I_{1,1}^T(\theta,\vec{q}\,)U_e + I_{1,2}(\theta,\vec{q}\,)f_\theta + E_1^T(\theta, \vec{q}\,) \dot{\vec{q}} .\nonumber \\
\dot{\psi}_3 &= \kappa_3(\vec{q}\,) \psi_3 + I_{3,1}^T(\theta,\vec{q}\,)U_e + I_{3,2}(\theta,\vec{q}\,)f_\theta + E_3^T(\theta, \vec{q}\,) \dot{\vec{q}}.
\end{align}
Note that because of the consideration of 3 adaptive parameters in the above Equations, $E_1$ and $E_3$ are both vectors of dimension 3.  Once again, $\dot{\vec{q}}$ and $f_\theta$ can be chosen arbitrarily.  Similar the structure of Equation \eqref{inveq} let
\begin{equation} \label{inveq2}
\begin{bmatrix}  \dot{q}_1 \\ f_\theta \\ \dot{q}_2 \\ \dot{q}_3  \end{bmatrix} = - A(\theta,\vec{q}\,)^{-1} 
    \begin{bmatrix}
    {\rm Real}(  I_{1,1}^T(\theta,\vec{q}\,)U_e ) \\ {\rm Imag}(  I_{1,1}^T(\theta,\vec{q}\,)U_e )  \\  {\rm Real}( I_{3,1}^T(\theta,\vec{q}\,)U_e)   \\ {\rm Imag}( I_{3,1}^T(\theta,\vec{q}\,)U_e)
    \end{bmatrix},
\end{equation} 
with
\begin{equation}
A(\theta,\vec{q}\,) = \begin{bmatrix}
 {\rm Real} ( E_{1,1}(\theta, \vec{q}\,))  & {\rm Real}(I_{1,2}(\theta,\vec{q}\,))  & {\rm Real}(E_{1,2}(\theta,\vec{q}\,)) & {\rm Real}(E_{1,3}(\theta,\vec{q}\,)) \\
  {\rm Imag} ( E_{1,1}(\theta, \vec{q}\,))  & {\rm Imag}(I_{1,2}(\theta,\vec{q}\,))  & {\rm Imag}(E_{1,2}(\theta,\vec{q}\,)) & {\rm Imag}(E_{1,3}(\theta,\vec{q}\,)) \\ 
   {\rm Real} ( E_{3,1}(\theta, \vec{q}\,))  & {\rm Real}(I_{3,2}(\theta,\vec{q}\,))  & {\rm Real}(E_{3,2}(\theta,\vec{q}\,)) & {\rm Real}(E_{3,3}(\theta,\vec{q}\,)) \\
      {\rm Imag} ( E_{3,1}(\theta, \vec{q}\,))  & {\rm Imag}(I_{3,2}(\theta,\vec{q}\,))  & {\rm Imag}(E_{3,2}(\theta,\vec{q}\,)) & {\rm Imag}(E_{3,3}(\theta,\vec{q}\,))
    \end{bmatrix},
\end{equation}
where $E_{1,j}$ and $E_{3,j}$ correspond to the $j^{\rm th}$ entries of $E_1$ and $E_3$, respectively.  With this choice for the parameter update values, Equation \eqref{psi1psi2} becomes $\dot{\psi}_j = \kappa_j(\vec{q}\,) \psi_j$ for $j = 1,2$.    Provided ${\rm Real}(\kappa_j(\vec{q}\,))<0$ for $j = 1,2$, both isostable coordinates tent to zero in the limit as $t$ approaches infinity and their dynamics can be ignored.  Or course, for \eqref{inveq} to be valid, $A(\theta,\vec{q}\,)^{-1}$ must exist for all $\vec{q}$ and $\theta$.  When $|| \vec{q} || = O(\epsilon)$, it is possible to compute  some of these terms directly by finding direct solutions to Equation \eqref{preg} and \eqref{irceq} to compute terms of the Floquet eigenfunctions and gradients of the Floquet coordinates, respectively, yielding the result
\begin{equation} \label{aaprox}
A(\theta,\vec{q}\;) \approx   \begin{bmatrix}
-1  & 0  &0 & 0\\
  0 & -q_1 ({\rm Imag}(\lambda_1)+\Delta \omega) & 0 & 0 \\ 
 0 & {\rm Real}(I_{3,2}(\theta,\vec{q}\,))  & -1 & 0 \\
    0  & {\rm Imag}(I_{3,2}(\theta,\vec{q}\,))  & 0 & -1
    \end{bmatrix}.
\end{equation}
Noting the diagonal structure of Equation \eqref{aaprox} invertibility is guaranteed when  $|| \vec{q}  \,|| = O(\epsilon)$.  

Taking $\dot{\vec{q}}$ and $f_\theta$ as mandated by Equation \eqref{inveq2} the adaptive reduction from \eqref{forcedadaptive} becomes
\begin{align} \label{largeadapt2}
\dot{\theta} &= \omega(q) + \omega(q) f_\theta, \nonumber \\
\dot{\vec{q}} &=  G_q(\vec{q},\theta,U_e), \nonumber \\
\dot{\psi}_j &= \kappa_j(\vec{q}\,) \psi_j + I^T_{j,1}(\theta,\vec{q}\,)U_e + I_{j,2}(\theta,\vec{q}\,) f_\theta + E_j^T(\theta,\vec{q}\,) \dot{\vec{q}\,}, \nonumber \\
j & = 5, \dots, \beta.
\end{align}
As compared to Equation \eqref{largeqadapt}, Equation \eqref{largeadapt2} simultaneously considers two nonlinear oscillatory modes but requires three total adaptive parameters.  As with the previous formulations, if $\kappa_j(\vec{q}\,)$ for any $j\geq 5$ is negative and large in magnitude for all allowable $\vec{q}$ it is generally possible to ignore the associated Floquet coordinate $\psi_j$ thereby yielding a reduced order model.

As a final note, it is straightforward to generalize the above strategy to consider more than two nonlinear modes.  However, this can become computationally prohibitive as each added nonlinear mode requires two additional adaptive parameters leading to an exponentially increasing amount of work required to compute the necessary family of periodic orbits and their associated Floquet eigenfunctions.

\section{Illustration of the Proposed Methodology} \label{ressec}

\subsection{Simple Pendulum}  \label{pendsec}
As a preliminary example meant to illustrate the implementation of the proposed strategy, consider the dynamics of a simple pendulum with viscous damping:
\begin{align} \label{pendeq}
\dot{x_1} &= x_2, \nonumber  \\
\dot{x_2} &= - \frac{g}{L} \sin(x_1) - \frac{b}{m L^2} x_2 + \frac{u(t)}{m L^2}.
\end{align}
Here, $x_1$ and $x_2$ correspond to the angular position, $\phi$ and velocity, $\dot{\phi}$, of a $m = 0.104$ Kg point mass suspended by a rigid, massless rod of length $L = 9.8$ m, $u$ is a torque input,  $g = 9.8 {\rm m}/{\rm s}^2$ is the acceleration due to gravity, and $b = 1\;{\rm Kg}\cdot {\rm m}^2/s$ is a viscous damping coefficient.   When $u = 0$ the pendulum has a stable fixed point at $x_1 = x_2 = 0$; eigenvalues of linearized fixed point are $\lambda_{1,2} = -0.050 \pm 0.999 i$.  We emphasize that the equations \eqref{pendeq} are already low dimensional and has a relatively simple nonlinearity. This example is intended to provide intuition about the implementation of the proposed strategy.

\begin{figure}[htb]
\begin{center}
\includegraphics[height=2.1 in]{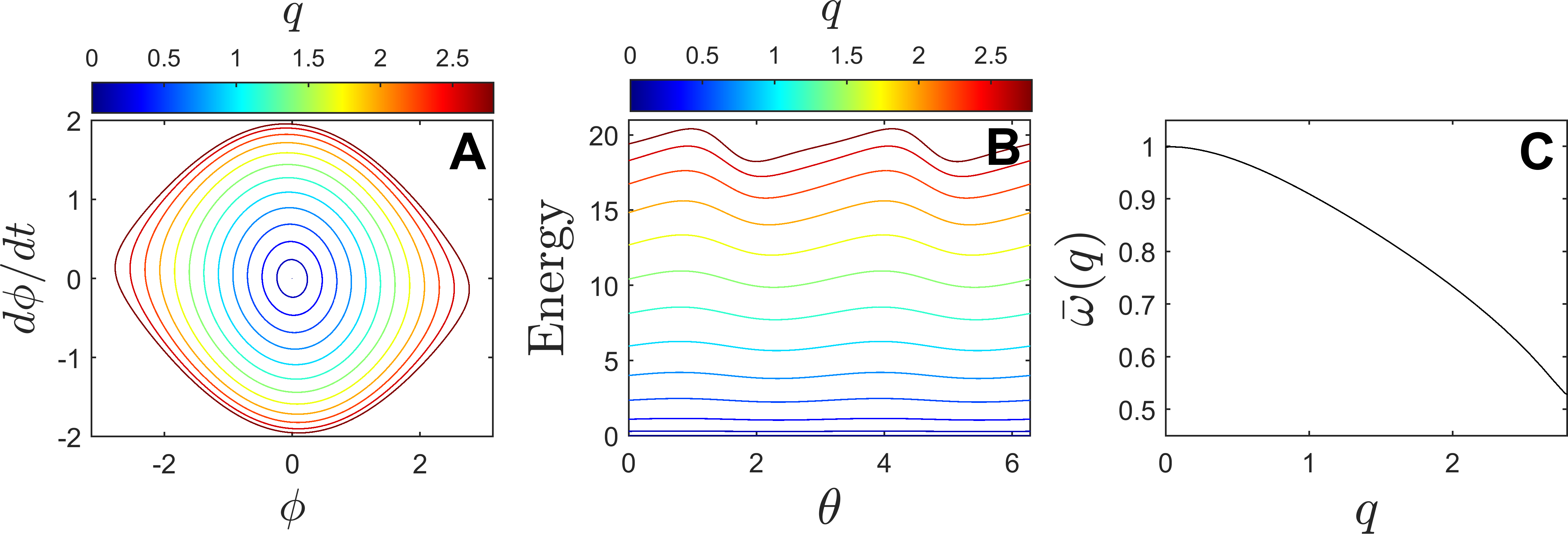}
\end{center}
\caption{The simple pendulum from \eqref{pendeq} is represented in the form \eqref{largeqadapt}.  Resulting forced periodic orbits are shown for different values of the amplitude-like variable $q$ in panel A. These orbits are reminiscent of (but not identical to) level sets of total energy as shown in panel B.  The effective natural frequency as a function of $q$ is shown in panel C.}
\label{penddetails}
\end{figure}

The proposed strategy as summarized in Section \ref{steplist} is applied to represent this model resulting in a phase-amplitude model of the form \eqref{largeqadapt}.  Note that for this specific example, this approach does not result in a reduction in dimension because the underlying model is already 2-dimensional. The resulting phase-amplitude model is comprised of different periodic orbits that emerge in response to external forcing, shown as a function of the amplitude-like parameter, $q$, in panel A of Figure \ref{penddetails}.   Despite their similar appearance, these orbits do not trace out level sets of total energy, i.e.,~$E = m g L(1-\cos(\phi)) + \frac{1}{2}m L^2 \dot{\phi}^2$ as shown in panel B.  The effective natural frequency computed according to \eqref{avfreq} is shown in panel C which is consistent with a lengthening period of oscillation for larger amplitude orbits.  In this example, it is not possible to continue the orbits beyond $\phi \approx \pm \pi$ which corresponds to a full revolution of the pendulum; near this point the resulting Floquet multipliers transition from complex-conjugate to real-valued rendering the iteration \eqref{nlinupdate} unusable beyond this point.

\begin{figure}[htb]
\begin{center}
\includegraphics[height=2.1 in]{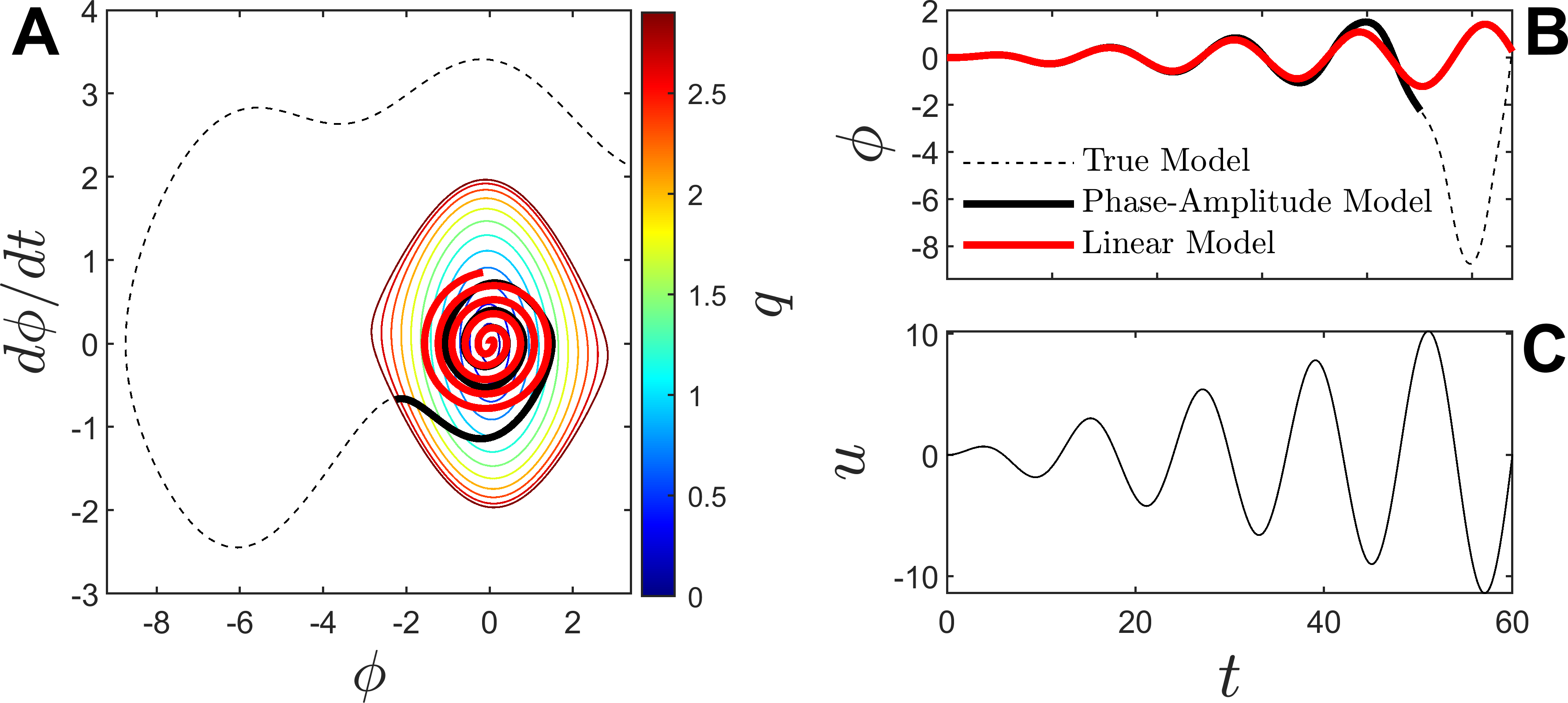}
\end{center}
\caption{Input of the form $u(t) =0.2 t \sin(\frac{2 \pi t}{12})$ is applied to the pendulum model from \eqref{pendeq} starting with an initial condition near steady state.  Black, red, and dashed lines show the output of the proposed phase-amplitude-based model, linearized model, and true model output in response to the input.  For reference, the colored lines show the forced periodic orbits that are used to define the phase-amplitude model.  Panel B gives a trace of the angular position for each model in response to the input plotted in Panel C.  In this example, the phase-amplitude model cannot be simulated beyond $t = 45$ because the state travels beyond the extent of the forced periodic orbits.  Nonetheless the timing and location of this traversal is nearly identical to the true model simulations.}
\label{pendresult}
\end{figure}

The accuracy of the phase-amplitude model is demonstrated in Figure \ref{pendresult}.  Panel A shows the output in response to the input $u(t) =0.2 t \sin(\frac{2 \pi t}{12})$ for the phase-amplitude model (black line) compared to the output of a model obtained through linearization about the fixed point (red line).  The true model output is shown as a dashed line.  The phase-amplitude model agrees perfectly with the full model output until the crossing of the outermost orbit, at which point, the model no longer displays predictable oscillations. Note that the phase-amplitude model cannot be used beyond this point because the state falls outside the family of forced periodic orbits; nonetheless, it pinpoints the exact moment that regular oscillations cease.  By contrast, the linear model matches the full model output for small amplitude oscillations but does not replicate the same sudden deviation from regular oscillations. Panel B of Figure \ref{pendresult} shows traces of $\phi$ and panel C shows the applied input.  These results are qualitatively similar when considering other inputs that drive the state past the regime that displays regular oscillations.

\FloatBarrier

\subsection{Coupled Population of Planar Oscillators}

Consider a model for a heterogeneous population of $N$ coupled planar oscillators
\begin{align} \label{isoclock}
\dot{x}_j &= \sigma x_j (\mu_j - r_j^2) - y_j (1+\rho_j (r_j^2 - \mu_j)) + \frac{K}{N}\sum_{j \neq i} x_j + u(t), \nonumber \\
\dot{y}_j &= \sigma y_j (\mu_j - r_j^2) + x_j (1+\rho_j ( r_j^2 - \mu_j)),
\end{align}
for $i = 1,\dots,N$ where $x$ and $y$ represent Cartesian coordinates, $r_j^2 = (x_j^2 + y_j^2)$, $N=10$,  $K = 1.2$ is the coupling strength, and $u(t)$ is an input common to each oscillator.  Additional model parameters are $\sigma = 0.1$, $\mu_j = -4 + 2j/9$, and $\rho_j = 0.4 - j/30$.    The individual elements in the model are similar to the radial isochron clock from \cite{winf01}; in the absence of coupling a stable limit cycle results from a Hopf bifurcation when $\mu_i>0$.  Here $\mu_i <0$ for all oscillators so that a stable fixed point results in this model at $x_j = y_j=0$ for all $j$ when $u(t)$ is held at 0.  

\begin{figure}[htb]
\begin{center}
\includegraphics[height=1.6 in]{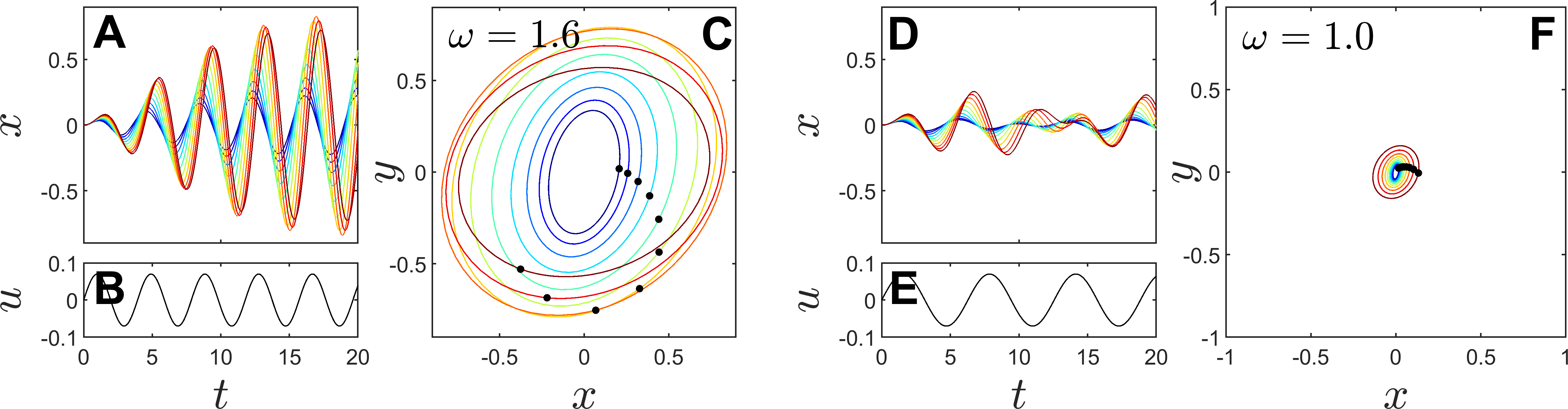}
\end{center}
\caption{Forced response for the coupled oscillator model from Equation \eqref{isoclock} in response sinusoidal input with frequency $\omega = 1.6$ (Panels A, B, and C) and $\omega = 1.0$ (Panels D, E, and F).  Panels A and D show traces of $x_i$ in response to the forcing in panels B and E, respectively, starting from steady state.  Panels C and F show the steady state response with different colors indicating the forced steady state orbit for each oscillator and dots give each oscillator's location at a snapshot in time.}
\label{sinusoidalforcing}
\end{figure}

Linearizing Equation \eqref{isoclock} about its fixed point gives 10 pairs of complex-conjugate eigenvalues.  The associated oscillatory modes have natural frequencies ranging between 1.23 and 2.57 rad/s.  Of particular interest is the eigenvalue pair $\lambda_{1,2} = -0.01 \pm 1.49i$ which has a slow decay rate relative to the next slowest decaying eigenvalue pair $\lambda_{3,4} = -0.25 \pm 1.23 i$.  This eigenvalue pair $\lambda_{1,2}$ causes a resonant peak for input frequencies near 1.5 rad/s; Figure \ref{sinusoidalforcing} shows the forced response (panels A and D) resulting from two different sinusoidal inputs (panels B and E).  Panels C and F give a representation for the steady state dynamics in response to periodic input $0.07 \sin(\omega t)$ taking $\omega = 1.6$ and $1.0$, respectively.  Orbits of different colors trace out the steady state solution for each oscillator and the dots provide a snapshot of each oscillator's relative position at a moment in time.  

Ultimately, the forced response peak is not well-captured by linear systems theory as illustrated in the results to follow.  Taking the output of the system \eqref{isoclock} to be $\bar{x} = \frac{1}{N}\sum_1^{N} x_i$, the system \eqref{isoclock} is linearized about its stable fixed point, and the resulting transfer function is used to predict the steady state response to inputs of various magnitude as shown in panel B.  Here, the amplitude is defined as $\max(\bar{x})-\min(\bar{x})$ over one period in steady state in response to sinusoidal forcing $u(t) = a \sin(\omega t)$.  As shown in panel B of Figure \ref{sinusoidalforcing}, the linearized model does not accurately capture the forced response near resonance.  Indeed, taking $a = 0.10$ predicts a peak height  that is more than 10 times larger than the true forced response (dashed lines).  Additionally, as the amplitude of forcing increases, the resonant peak shifts towards faster frequencies, a feature that cannot be captured by any linear model.  As an alternative approach, the reduced order modeling strategy described in Section \ref{steplist} is applied to the model \eqref{isoclock} to obtain a model of the form \eqref{largeqadapt}.  The resulting model truncates all Floquet coordinates $\psi_3,\dots,\psi_{N}$ associated with the faster decaying modes; as such, the resulting nonlinear model is 2-dimensional.  Panel A of Figure \ref{sinusoidalforcing} shows the steady state response of the forced reduced order model, providing a much more accurate match than the linear model.  For the nonlinear model, the curves are obtained by identifying periodic solutions of \eqref{largeqadapt} in response to the indicated input.

\begin{figure}[htb]
\begin{center}
\includegraphics[height=1.8 in]{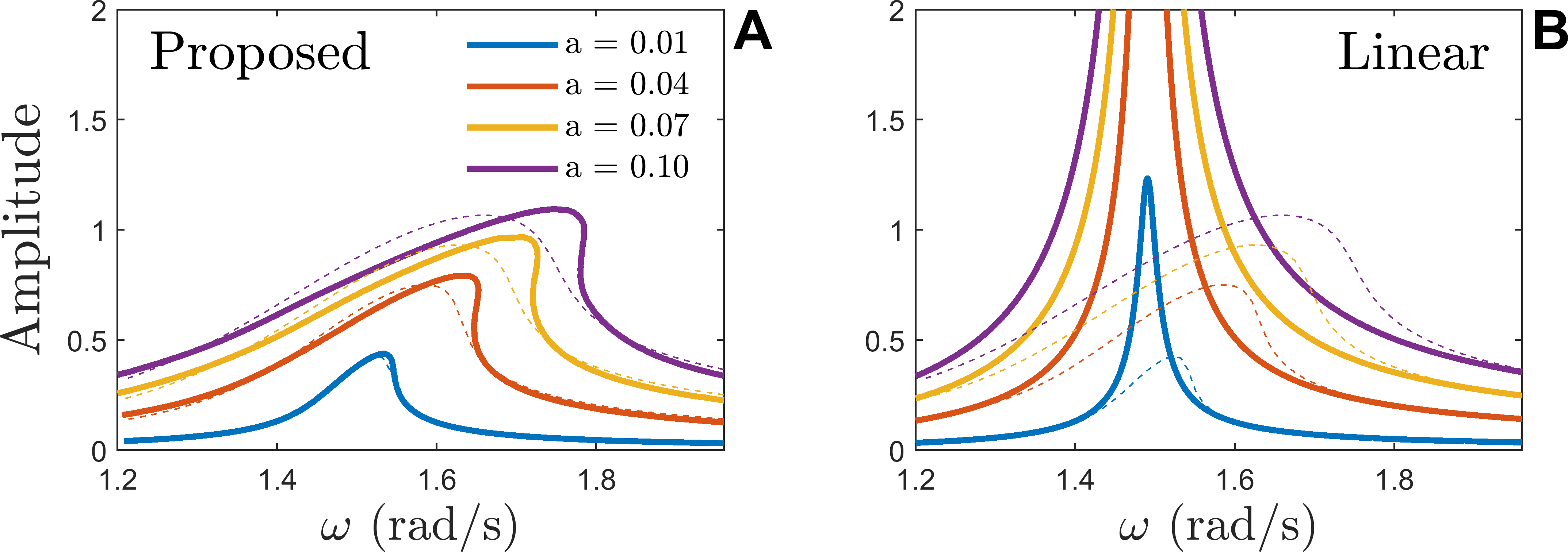}
\end{center}
\caption{ The steady state response to sinusoidal forcing $u(t) = a \sin(\omega t)$ is considered for the coupled oscillator model from \eqref{isoclock}.  Amplitude is defined as the difference between the maximum and minimum value of $\bar{x} = \frac{1}{N}\sum_1^{N} x_i$ in steady state.  Accuracy of the proposed reduced order model of the form \eqref{largeqadapt} is compared to a linear model obtained by linearizing \eqref{isoclock} about its fixed point with results shown in panels A and B, respectively.  The reduced order model truncates all Floquet coordinates $\psi_3,\dots,\psi_{N}$ so that there are only 2 dimensions.    The dashed lines show the true steady state response of the full model obtained from simulations of Equation \eqref{isoclock}.}
\label{forcedresults}
\end{figure}

\begin{figure}[htb]
\begin{center}
\includegraphics[height=2.0 in]{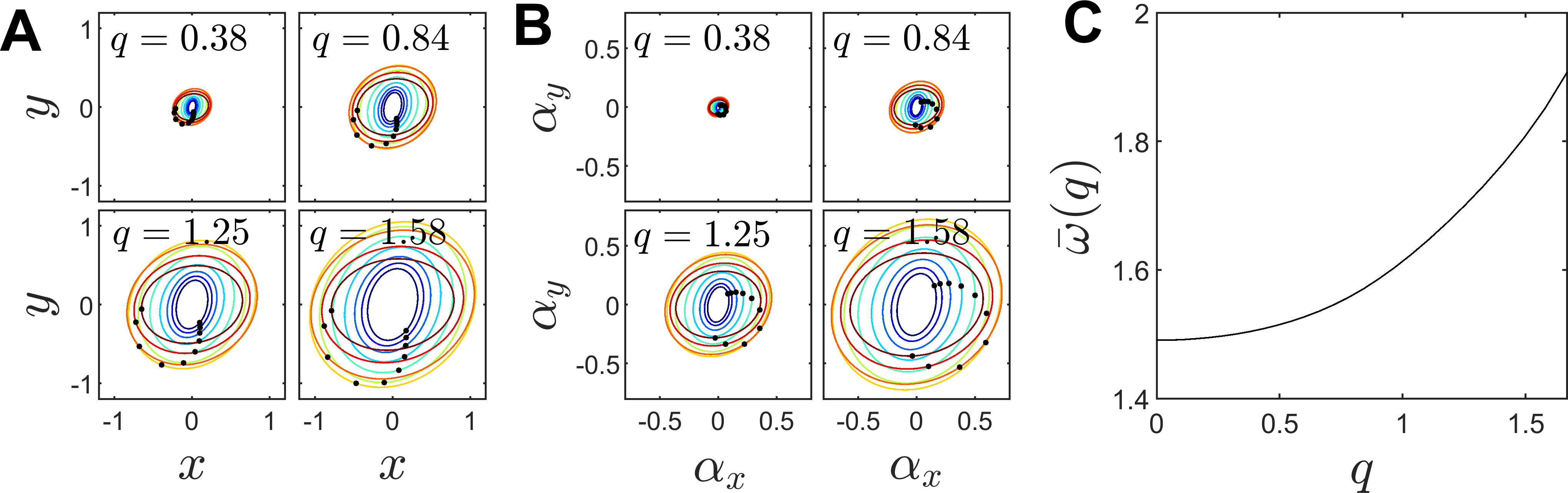}
\end{center}
\caption{ The coupled oscillator model from Equation \eqref{isoclock} is represented in the form \eqref{largeqadapt}.  Forced periodic orbits for different values of the amplitude-like variable $q$ are shown in panel A with different colors corresponding to the orbits of different oscillators.  The applied input $\alpha(t,q)$ is represented in panel B for different values of $q$.  The different colored traces correspond to the input applied to each oscillator over a single oscillation.  Black dots provide snapshots of the state and applied input at a single moment in time for the periodic orbit.  Panel C illustrates that the effective natural frequency computed according to \eqref{avfreq} increases with oscillation amplitude.}
\label{hopforbits}
\end{figure}

Figure \ref{hopforbits} provides additional information about the resulting reduced order model of the form \eqref{largeqadapt}.  For different values of $q$, panel A shows traces of the orbit of each oscillator for the forced periodic orbit $x^\gamma_q$ associated with the input $\alpha(q,t)$.  Panel B shows the associated input $\alpha(q,t)$.  In panel A, different colors correspond to the the forced periodic orbit of a given oscillator.  In panel B, different colors correspond to the periodic forcing applied to given oscillator, i.e.,~with $\alpha_x$ and $\alpha_y$ representing the input applied to the $x$ and $y$ coordinate, respectively, that yields the periodic solution.  Note that while the input $u(t)$ from Equation \eqref{isoclock} only appears in the $x$-coordinate dynamics of each oscillator, the input $\alpha(q,t)$ is applied to both the $x$ and $y$ coordinates for each oscillator.  Black dots in panels A and B provide snapshots of the state and applied inputs at a given moment in time.  The periodic solutions for each oscillator in panel A are similar, but not identical, to the orbits that emerges in response to sinusoidal input as shown in Panel C of Figure \ref{hopforbits}.  Large amplitude oscillations are accurately captured by this single mode for inputs that are near the resonant frequency.  Panel C of Figure \ref{hopforbits} gives the effective unforced oscillation frequency for the reduced order model of the form \eqref{largeqadapt} for different values of $q$ computed according to Equation \eqref{avfreq}.  For this model the effective natural frequency grows as the amplitude of the oscillation increases.

\FloatBarrier

\subsection{Power System Model}
Finally, we consider a mathematical model of a power system comprised of ordinary differential equations on synchronous generators, controllers and other dynamic devices and algebraic equations on the power network. The classical $m$-generator system model  (\ref{eq_powersystem}) is used here, in which the $i$-th generator is represented by its swing equations including two first-order differential equations respectively on its rotor angle $\delta_i$ (in radians) and rotor speed $\omega_i$ (in radian per second), and an algebraic equation couples its electric power output $P_{e,i}$ with the rotor angles of all other networked generators.

\begin{equation}
	\begin{array}{lll}
		\label{eq_powersystem}
		\dot{\omega_i} =\frac{\omega_{0}}{2H_{i}}\left( P_{m,i}-D_i\frac{\omega_{i}-\omega_{0}}{\omega_{0}}- P_{e,i} \right), \\
		\dot{\delta_i}  =\omega_{i}-\omega_{0}, \\
  	P_{e,i}=E_i^2G_{ii}+\displaystyle\sum_{j=1, j\neq i}^m E_iE_j(G_{ij}\cos(\delta_i-\delta_j)+B_{ij}\sin(\delta_i-\delta_j)), &1\leq i\leq m.
	\end{array}
\end{equation}

\noindent Above, $\omega_{0}$ is the system's synchronous speed. Other parameters on the $i$-th generator include its mechanical power input $P_{m,i}$ from the turbine, which is considered constant, inertial time constant $H_{i}$, damping coefficient $D_i$, and electromotive force $E_i$ treated as a constant under the excitation control. All branches and loads of the power network are equivalenced by constant admittances such as $G_{ij}+jB_{ij}$ between two generators and conductances to the ground such as $G_{ii}$. Further details about the model and its parameters are given in \cite{anderson}. Here we consider the 3-generator power system model with $m=3$ in  \cite{anderson}, which is the so-called IEEE 3-generator 9-bus test system. 

\begin{figure}[htb]
\begin{center}
\includegraphics[height=3.2 in]{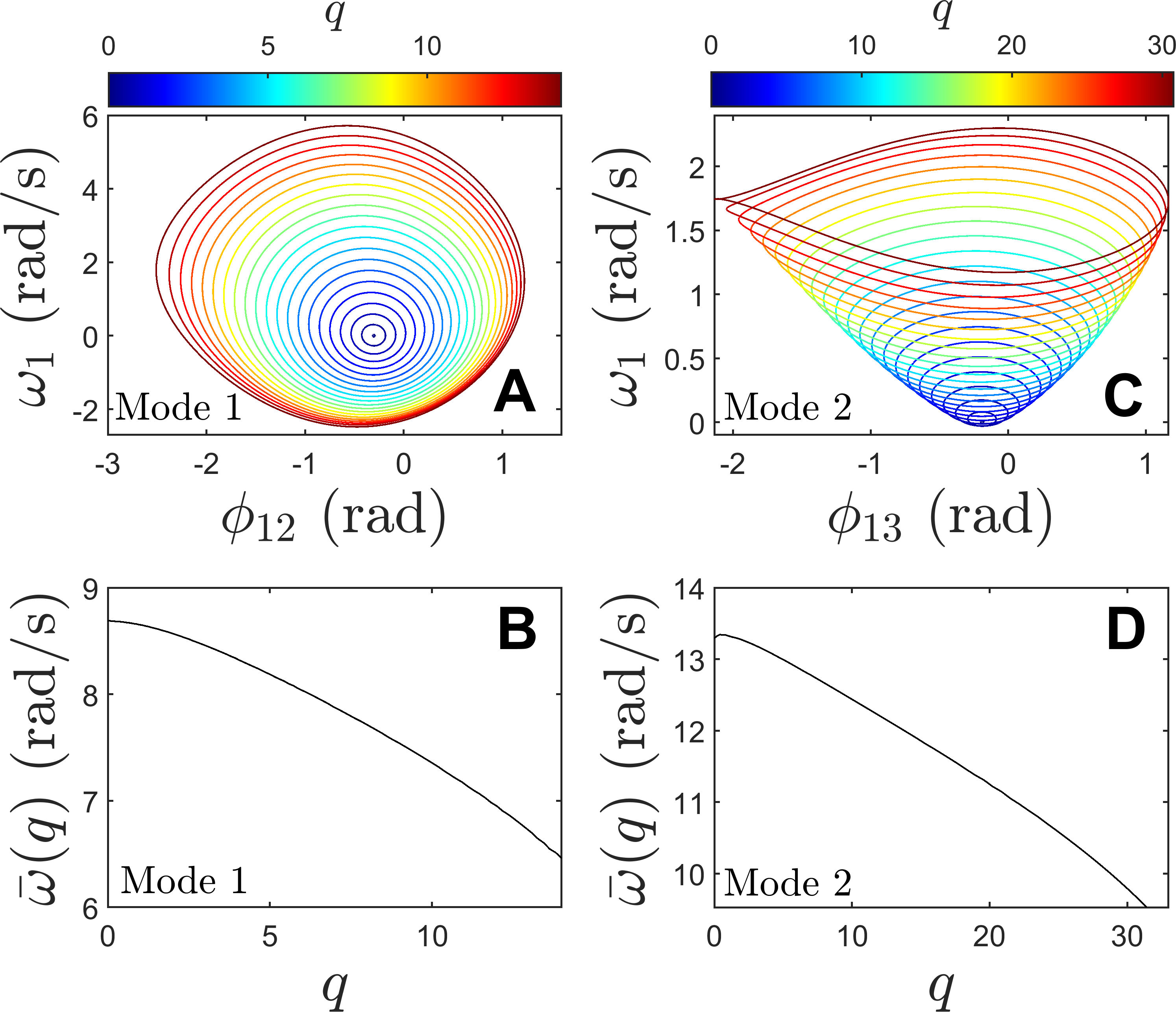}
\end{center}
\caption{Two-dimensional projections of the forced periodic orbits associated with slower and faster nonlinear oscillations (mode 1 and mode 2, respectively) are shown in panels A and C.  For small values of the amplitude-like parameter $q$, these are nearly identical to linear oscillations governed by their respective eigenvalues and eigenvectors.  Larger amplitude modes associated with larger values of $q$ cannot be accurately represented as linear oscillations.  Panels B and D show the effective frequencies computed according to Equation \eqref{avfreq} for each oscillation mode as a function of $q$.  For each mode, the effective frequency decreases as the oscillation magnitude increases.}
\label{psmodes}
\end{figure}

As a preliminary step in the analysis, the dynamics of Equation \eqref{eq_powersystem} are considered in reference to angle $\delta_1$ defining $\phi_{12} = \delta_1 - \delta_2$ and $\phi_{13} = \delta_1 - \delta_2$. The dynamics of the phase differences are given by
\begin{equation}
\dot{\phi}_{1j} = \omega_1 - \omega_j,
\end{equation}
for $j = 2,3$.   The rotor speed equations can be written as a function of $\phi_{12}$ and $\phi_{13}$,  yielding a 5-dimensional model with a stable fixed point at $\begin{bmatrix} \phi_{12} & \phi_{13} & \omega_1 & \omega_2 & \omega_3 \end{bmatrix} = \begin{bmatrix}-0.30 & -0.19 & 0 & 0 & 0 \end{bmatrix}$.  In steady state, the rotors are phase cohesive with identical frequencies and slight differences between their phases.   The eigenvalues associated with this fixed point are  $\lambda_{1,2} = -0.25 \pm 8.69i$, $\lambda_{3,4} = -0.25 \pm 13.36 i$ and $\lambda_5 = -0.5$.  The first four eigenvalues correspond to oscillatory eigenmodes with frequencies of 1.38 and 2.12 Hz.  The slow and fast oscillatory mode will be referred to as mode 1 and mode 2, respectively.  Individually, these modes are used to obtain two different reduced order models of the form \eqref{largeqadapt} following the proposed strategy summarized in Section \ref{steplist}.  The first (resp.,~second) model can accommodate large amplitude, nonlinear mode 1 (resp.,~mode 2) oscillations.  Two-dimensional projections of the associated periodic orbits for different values of the amplitude-like parameter $q$ are shown in panel A (resp.,~C) of Figure \ref{psmodes} with effective frequencies shown in Panel B (resp.,~D) computed according to Equation \eqref{avfreq}.   It is not possible to continue the periodic orbits beyond the limits shown in panels A and C of Figure \ref{psmodes}; beyond this point their Floquet multipliers transition from complex-conjugate to real-valued so that the iteration from Equation \eqref{nlinupdate} (which is implemented as part of step 7 of the procedure from Section \ref{steplist}) cannot be continued.  Note that much like for the simple pendulum example from Section \ref{pendsec}, this limit coincides with a qualitative change in the collective behavior of the model \eqref{eq_powersystem} transitioning to regions of phase space that do not yield oscillatory dynamics.

\begin{figure}[htb]
\begin{center}
\includegraphics[height=2.4 in]{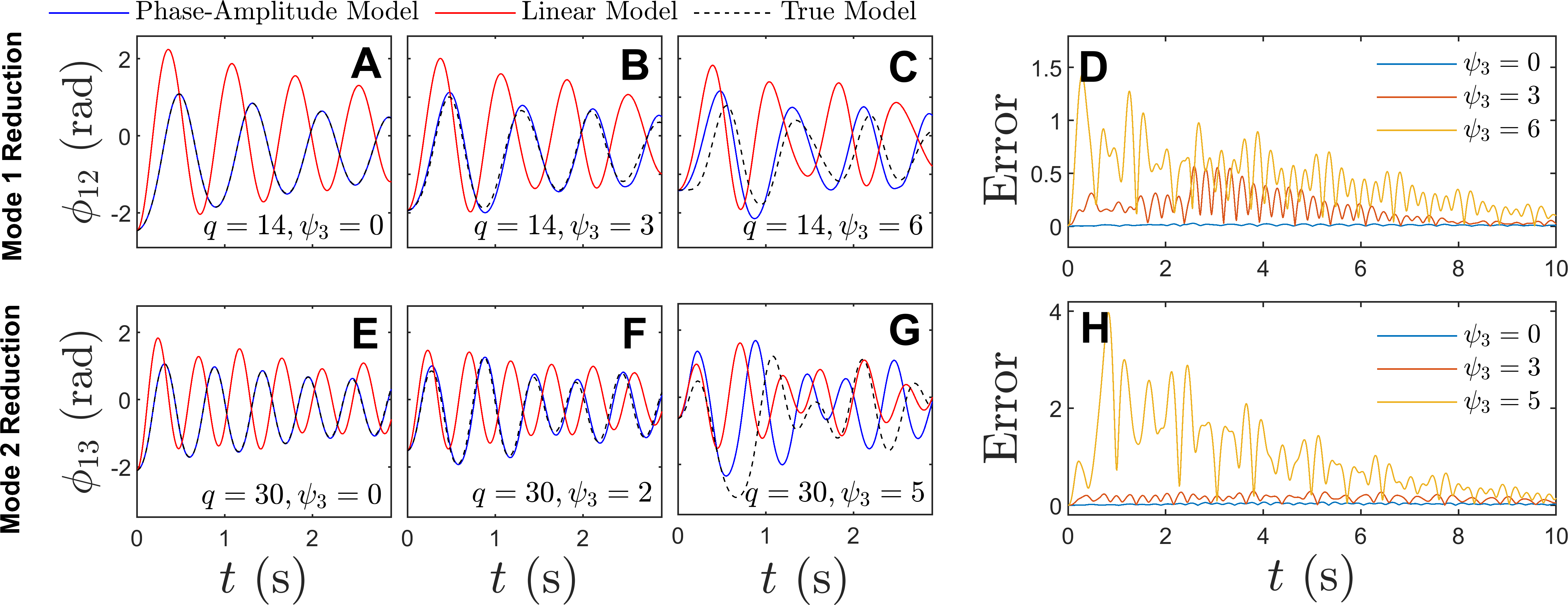}
\end{center}
\caption{ Panels A, B, and C (resp.,~E, F, and G) show dominant (nonlinear) oscillations from mode 1 (resp.,~mode 2)  with varying contributions from the non-dominant mode mode.  The magnitude of $q$ correlates with the amplitude of the dominant nonlinear mode and the magnitude of $\psi_3$ correlates with amplitude of the non-dominant mode.  Initial conditions indicated on each panel are chosen to highlight the effect of increasing the non-dominant mode, i.e.,~by taking $\psi_3$ to be larger.  Simulation results from a linear model are also provided, obtained by linearizing about the stable fixed point, to demonstrate that this state is beyond the regime for which linearization is accurate.  Panels D and H show the two norm of the error between the phase differences for the the proposed phase-amplitude model \eqref{largeqadapt} and the true model \eqref{eq_powersystem} for different initial values of $\psi_3$.  When the non-dominant mode is small, the prediction accuracy is nearly flawless.  As the non-dominant mode grows, the accuracy suffers.}
\label{psresult}
\end{figure}

Each of the models described in Figure \ref{psmodes} can accommodate large oscillations associated with a single nonlinear mode.  The variables $\theta$ and $q$ in Equation \eqref{largeqadapt} capture oscillations of the dominant mode.  Oscillations associated with the other mode are captured using complex-conjugate Floquet coordinates $\psi_3$ and $\psi_4$ which are assumed to be of small magnitude.  Figure \ref{psresult} illustrates the accuracy of the resulting reduced order models in relation to this assumption.  Panels A, B, and C show nonlinear mode 1 oscillations with varying contributions from mode 2 (quantified by the value of $\psi_3$ at $t = 0$).  When $\psi_3 = 0$, the reduced order simulation of the model \eqref{largeqadapt} (blue line) gives results that are indistinguishable to those from the full order simulation of \eqref{eq_powersystem} (dashed line).  As the contribution from mode 2 (and hence the magnitude of the initial value of $\psi_3$) increases, the error between the full order and reduced order models increases.  For comparison, simulation results are also provided for a model obtained by linearizing the dynamics about the stable fixed point (red lines); the oscillations considered here are clearly beyond the regime for which linearization provides an accurate representation for the system dynamics.  Panel D shows the two-norm of the error associated with the variables $\phi_{12}$ and $\phi_{23}$ for the phase-amplitude model for different initial values of $\psi_3$.  Panels E, F, and G, provide analogous results to panels A, B, and C, except when considering nonlinear mode 2 oscillations and requiring the mode 1 oscillations to be small.   Likewise, panel H shows the resulting error between the phase-amplitude and full order models for differing contributions from the non-dominant mode.

\begin{figure}[htb]
\begin{center}
\includegraphics[height=2.8 in]{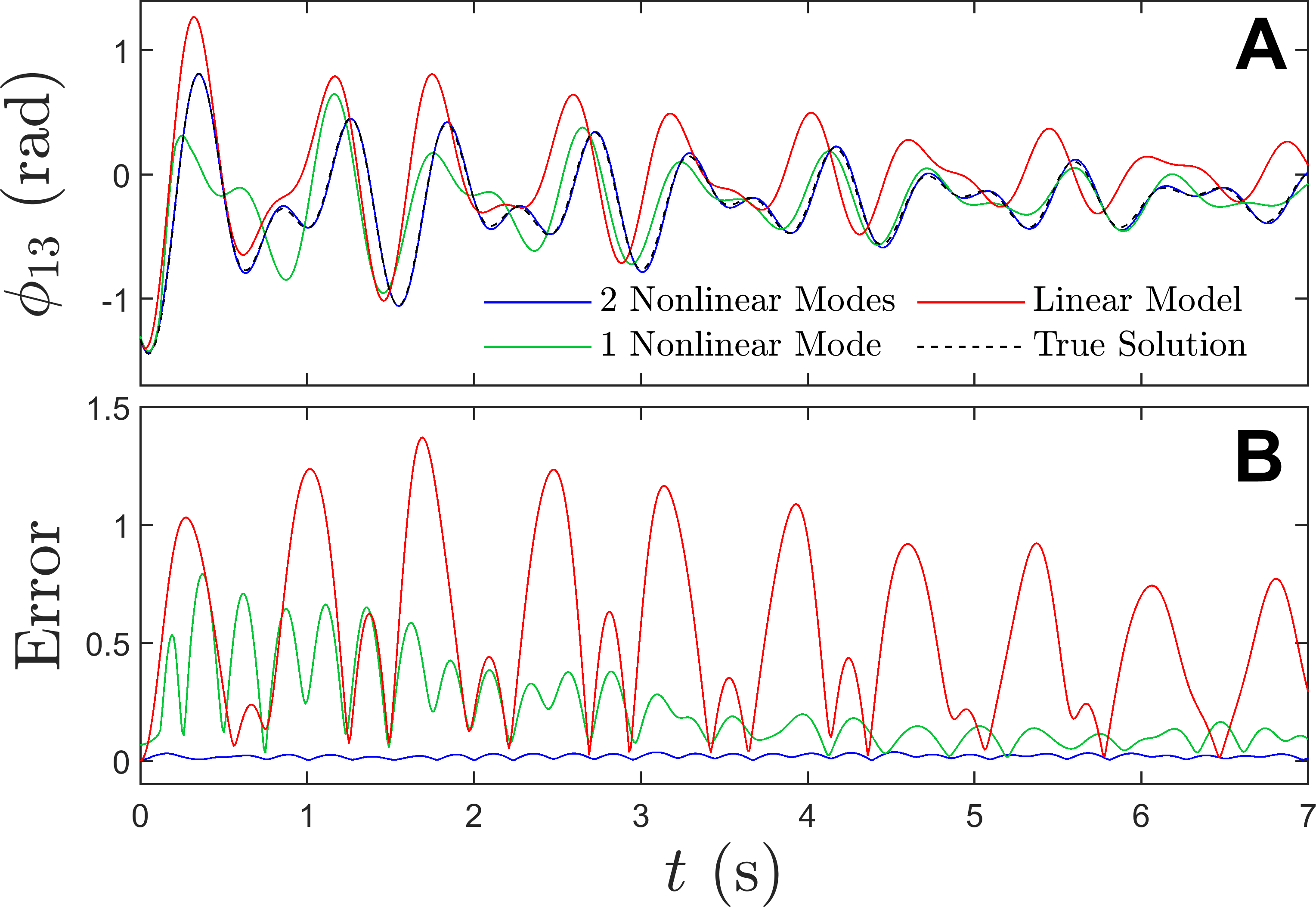}
\end{center}
\caption{Reduced order simulations of \eqref{eq_powersystem} considering two nonlinear oscillatory modes.  For an initial condition asssociated with moderate contributions from both mode 1 and mode 2, a trace of $\phi_{13}$ is shown in panel A for the two-nonlinear-mode model \eqref{largeadapt2} and compared to outputs from the true model (blue and dashed line, respectively).  Comparable simulations of a one-nonlinear-mode model of the form \eqref{largeqadapt} and a model obtained through local linearization of \eqref{eq_powersystem} about its fixed point are shown with green and red lines, respectively.  Panel E shows the two norm of the error between the phase differences for each of the models considered.}
\label{twomodes}
\end{figure}

The method from Section \ref{multimodes} is also applied to simultaneously consider two nonlinear oscillation modes. Here, mode 1 (i.e.,~the slow oscillation mode stemming from the eigenvalues $\lambda_{1,2} = -0.25 \pm 8.69i$) is used to define periodic orbits $x^\gamma_{q_1,0,0}(t)$ with associated inputs $\alpha(q_1,0,0)$ as described in Section \ref{multimodes}.  A 3-parameter family of orbits is defined iteratively using Equations \eqref{r1eq}-\eqref{finalalpha}.  For the implementation of this strategy, $g_3$ corresponds to the Floquet eigenfunction associated with mode 2 (i.e.,~the fast oscillation mode stemming from the eigenvalues $\lambda_{3,4} = -0.25 \pm 13.36 i$).  This information is used to define a model of the form \eqref{largeadapt2} that contains one phase coordinate and 3 amplitude-like coordinates $q_1,q_2$, and $q_3$, and no additional Floquet coordinates for a total of 4 dimensions.  This model can consider oscillations with large contributions from both mode 1 and mode 2.  By contrast, models of the form \eqref{largeqadapt} that were used to obtain results in Figure \ref{psresult} can only accurately consider nonlinear oscillations with large magnitude contributions from either mode 1 or mode 2, but not both.  Results in Figure \ref{twomodes} highlight this distinction.  The two-nonlinear-mode model is simulated using an initial condition $\begin{bmatrix} \theta & p_1 & p_2 & p_3 \end{bmatrix} = \begin{bmatrix} 0 & 5.1 & 2.1 & -2.1 \end{bmatrix}$ that corresponds to a state that yields oscillations with moderate contributions from both mode 1 and mode 2.  In panel A of Figure \ref{twomodes}, output from the two-nonlinear-mode model (blue line) is compared to the output from the true model (dashed line) with results that are nearly indistinguishable.  A comparable initial condition is used in a simulation of the one-nonlinear-mode model (i.e.,  the same model from panels A-D of Figure \ref{psresult} with results shown in green; this model is not able to accurately accommodate the contribution from both modes.  For reference, simulation results when using a model obtained from linearization of the stable fixed point are shown in red which performs worse than the other two models.  Panel E shows the two norm of the error between the phase differences for each of the models considered. The results presented here are qualitatively similar when considering other initial conditions.

\FloatBarrier

\section{Discussion and Conclusion} \label{concsec}

This work leverages recently developed adaptive phase-amplitude reduction techniques \cite{wils21adapt} to yield a general approach for reduced order characterization of systems with highly nonlinear oscillations.  By carefully defining a family of periodic orbits associated with a particular mode of oscillation, a reduced order model of the form \eqref{largeqadapt} can be obtained that accurately replicates large amplitude nonlinear oscillations that emerge in response to arbitrary external forcing;  in the limit that the oscillation amplitude is small, the resulting modes of the proposed strategy become functionally identical to linear modes.  This approach can also be used to consider interactions between multiple nonlinear modes as discussed in Section \ref{multimodes}.   In the examples considered in Section \ref{ressec}, the proposed reduced order modeling strategy provides substantially improved results as compared to non-reduced order models obtained using local linearization.

It would be of general interest to more carefully investigate the relationship between the proposed approach and other nonlinear model representation techniques.  The proposed approach shares similarities with the notion of spectral submanifolds described in detail in \cite{hall16}.  In the limit as the state relaxes to the fixed point, both methods yield models that behave similarly to those obtained from linear approximation.  Furthermore, both techniques can exploit differences in spectral gaps between the slowest decaying eigenvalues to obtain reduced order models.  In contrast, however, while spectral submanifolds are invariant under the flow, the family of periodic orbits used to construct the reduced order models (and ultimately used to define the nonlinear modes) is not an invariant set.  Additionally, spectral submanifolds are defined in relation to the \emph{unforced} dynamics of the underlying system whereas external forcing is inherent to the implementation of the proposed approach.  As far as other approaches,  while the proposed strategy does not explicitly consider the notion of isostable coordinates associated with a fixed point \cite{maur13}, \cite{maur16} (i.e.~level sets of the slowest decaying modes of the Koopman operator) there is also a possible connection with the amplitude coordinates used in the proposed strategy from Equation \eqref{largeqadapt} as the dynamics are identical in the limit that the state approaches the fixed point.

While the results of this work are promising, there are a number of limitations left to address.  It is necessary to know the underlying equations in order to implement the proposed approach, for example, in order to obtain solutions of \eqref{adjiso} and \eqref{adjfloq} in the computation of the Floquet eigenfunctions and the gradients of the Floquet coordinates.  Additional modifications would be necessary to implement this approach in a data-driven setting.  Additionally, in principle this approach can be used to consider the interactions between an arbitrary number of nonlinear modes, the computational expense associated with computing the required periodic orbits grows quickly with the number of modes considered.  Indeed, the iterations associated with the two-mode reduction from Equations \eqref{r1eq}-\eqref{finalalpha} must be performed in three dimensions instead of the single dimension required when considering only one mode.  When considering a system with many nonlinear modes, it would likely be necessary to consider multiple separate models that characterize the interactions between smaller subsets of nonlinear modes of interest.  Finally, in the examples considered in this work, we were not able to consider orbits that extended beyond the basin of attraction of the fixed point.  Particularly for both the simple pendulum \eqref{pendeq} and the power system model \eqref{eq_powersystem}, as the states approach the boundary of the basin of attraction of the fixed point, the Floquet multipliers associated with the periodic orbits become real-valued precluding the continuation of the iteration used for defining the family of periodic orbits.  It would be useful to develop a workaround for this issue in order to extend the applicability of this approach.

This material is based upon the work supported by the National Science Foundation (NSF) under Grant No.~CMMI-2140527.

\begin{appendices}

\section{Computation of the Terms Comprising the Phase and Phase-Amplitude Reduced Order Equations} \label{apx0}
\renewcommand{\thetable}{A\arabic{table}}  
\renewcommand{\thefigure}{A\arabic{figure}} 
\renewcommand{\theequation}{A\arabic{equation}} 
\setcounter{equation}{0}
\setcounter{figure}{0}

The gradient of the phase from the phase reduction \eqref{predeq} must generally be computed numerically.  This problem has been studied widely \cite{brow04}, \cite{erme10}, \cite{hopp97}; holding $q$ constant,  solutions of $Z(\theta,q)$ can be obtained by finding periodic solutions of the adjoint equation
\begin{equation} \label{adjeq}
\dot{Z} = -J^T Z,
\end{equation}
where $J$ is the Jacobian evaluated at $y_q^\gamma$.  Equation \eqref{adjeq} always has a single non-decaying Floquet exponent. As such periodic solutions must be normalized so that $\frac{2 \pi}{T(q)} = F^T(y^\gamma_q(t),q) Z(\theta(t),q)$.  As detailed in \cite{wils16isos}, $I_j(\theta,q)$ as defined in the phase-amplitude transformed equations \eqref{phaseamp}  (i.e.,~the gradient of the $j^{\rm th}$ Floquet coordinate with respect to the state) can also be obtained by finding periodic solutions of 
\begin{equation} \label{adjiso}
\dot{I}_j = -(J^T - \kappa_j {\rm Id}) I_j,
\end{equation}
where $\kappa_j$ is the Floquet exponent associated with $\psi_j$ and ${\rm Id}$ is the identity matrix of appropriate size.  Likewise, the Floquet eigenfunctions $g_j(\theta,q)$ can be obtained by finding periodic solutions of 
\begin{equation} \label{adjfloq}
\dot{g}_j = (J-\kappa_j {\rm Id}) g_j.
\end{equation}
As discussed in \cite{wils20highacc}, solutions of Equations \eqref{adjiso} and \eqref{adjfloq} follow the relationships
\begin{equation} \label{gconst}
g_k^T(\theta,q) I_j(\theta,q) =  \begin{cases}  1   &  \text{if } k = j,    \\
   0 ,  & \text{otherwise}.  \end{cases}
\end{equation}
Additionally,
\begin{equation} \label{normeq}
I_j^T(\theta,q) \frac{\partial y^\gamma_q}{\partial \theta} = 0 \text{ for all } j.
\end{equation}
When considering the adaptive phase-amplitude reduced order equations from Equation \eqref{phaseampred}, reference \cite{wils21adapt} established direct relationships between  $Z(\theta,q)$ (resp.,~$I_j(\theta,q)$) and the term $D(\theta,q)$ (resp.,~$E_j(\theta,q)$).  Specifically, letting $\frac{\partial y^\gamma}{\partial q_j} |_{\theta_0,q} \equiv \lim_{a \rightarrow 0}  (y^\gamma_{q + e_j a}(\theta_0) - y^\gamma_{q}(\theta_0))/a$ where $e_j$ is the $j^{\rm th}$ component of the standard unit basis, one can show that
\begin{equation}
e_k^T D(\theta,q) = -Z^T(\theta,q) \frac{\partial y^\gamma_q}{\partial q_k}, 
\end{equation}
and
\begin{equation} \label{eeq}
e_k^T E_j(\theta,q) = -I_j^T(\theta,q) \frac{\partial y^\gamma_q}{\partial q_k}.
\end{equation}

\section{First Order Perturbations of Simple Eigenvalues and Corresponding Eigenvectors} \label{apxb}

\renewcommand{\thetable}{B\arabic{table}}  
\renewcommand{\thefigure}{B\arabic{figure}} 
\renewcommand{\theequation}{B\arabic{equation}} 
\setcounter{equation}{0}
\setcounter{figure}{0}

Let $A \in \mathbb{R}^{N \times N}$ have a simple (i.e.,~unique) eigenvalue $\lambda$ with corresponding left and right eigenvectors $w$ and $v$, respectively, normalized so that $w^* v = 1$, $v^*v = 1$, and ${\rm arg}(e_j^T v_1) =-\pi$.  Here, $e_j$ is the $j^{\rm th}$ element of the standard unit basis, ${\rm arg}(\cdot)$ is the argument of the complex number, $j$ can be chosen arbitrarily, and $^*$ denotes the conjugate transpose.  By definition, $\lambda$ and $v$ solve
\begin{equation} \label{eval}
f(A,v,\lambda) = 0 = Av - \lambda v.
\end{equation}
The goal is to characterize the change in the eigenvalue $\lambda + d\lambda$ and eigenvector $v + dv$ that result when the matrix $A$ is shifted incrementally to $A + dA$.  Taking the total differential of Equation \eqref{eval} yields
\begin{equation} \label{mydiff}
df = 0 = dA v - d\lambda v + Adv-\lambda dv,
\end{equation}
which must be satisfied for the perturbed eigenvalue/eigenvector pair $(\lambda + d\lambda, v + dv)$.  Multiplying on the left by $w^*$ and rearranging Equation \eqref{mydiff} yields
\begin{equation}
w^* d\lambda v = w^* dA v + w^* A dv - w^* \lambda dv.
\end{equation}
Using the fact that $w^* v = 1$ and $w^* A = w^* \lambda$, the above equation simplifies to 
\begin{equation}  \label{dleq}
d \lambda = w^* dA v.
\end{equation}
Again considering Equation \eqref{mydiff}, the perturbation in the eigenvector can be obtained by solving
\begin{equation} \label{dveq}
(  A - \lambda {\rm Id} ) dv =   d \lambda v - dA v,
\end{equation}
where ${\rm Id}$ is an appropriately sized identity matrix.  Equation \eqref{dveq} above is obtained through manipulation of \eqref{mydiff}.  Noting that $v$ is in the null space of $A-\lambda {\rm Id}$, Equation \eqref{dveq} only has solutions if $w$ is orthogonal to $d \lambda v - dA v$.  One can verify this is the case directly:
\begin{align}
w^*(d \lambda v - dA v) &= w^*v (w^* dA v)    - w^* dA v \nonumber \\
&= 0,
\end{align}
where the right hand side of the first line is obtained by substituting Equation \eqref{dleq} and the second line is obtained by changing the order of multiplication noticing that $w^* dA v \in \mathbb{C}$  and recalling that $w^* v = 1$.  Thus all solutions of Equation \eqref{dveq} are given by
\begin{equation}
dv = (A-\lambda {\rm Id})^\dagger   (d\lambda v - dA v) + v \alpha,
\end{equation}
for any $\alpha \in\mathbb{C}$ where $^\dagger$ denotes the Moore-Penrose pseudoinverse.  Here, $\alpha$ must be chosen appropriately so that $v + dv$ satisfies the required normalization.

\end{appendices}

\end{document}